\documentclass[final,onefignum,onetabnum,letterpaper]{siamart250211}

\usepackage{soul}
\usepackage{lipsum}
\usepackage{amsfonts}
\usepackage{epstopdf}
\usepackage{epsfig}
\ifpdf
  \DeclareGraphicsExtensions{.eps,.pdf,.png,.jpg}
\else
  \DeclareGraphicsExtensions{.eps}
\fi

\usepackage{datetime}
\newdateformat{monthyeardate}{%
	\monthname[\THEMONTH] \THEDAY, \THEYEAR}

\usepackage{academicons}
\usepackage{xcolor}

\usepackage{amsmath}
\usepackage{amssymb}
\usepackage{commath}
\usepackage{mathtools}
\usepackage{bbm}
\usepackage{bm}
\usepackage{upgreek}
\usepackage{float}
\allowdisplaybreaks
\usepackage{import}
\usepackage{color}
\usepackage{graphicx}
\usepackage[small]{caption}
\usepackage{subcaption}
\usepackage{tabularray}
\UseTblrLibrary{amsmath}

\usepackage{relsize}
\usepackage{adjustbox}
\usepackage{algorithm}
\usepackage{algpseudocode}
\usepackage{booktabs}
\usepackage{tikz}
\usepackage{pifont}
\usetikzlibrary{bayesnet} 

\usepackage{enumitem}
\setlist[enumerate]{leftmargin=.5in}
\setlist[itemize]{leftmargin=.5in}

\newsiamremark{remark}{Remark}
\crefname{hypothesis}{Hypothesis}{Hypotheses}

\headers{Entropy-stable learned Hyperbolic Conservation Laws}{Lizuo Liu, Lu Zhang and Anne Gelb}

\title{Neural Entropy-stable conservative flux form neural networks for learning hyperbolic conservation laws
\thanks{
\monthyeardate\today 
}}

\author{
Lizuo Liu\thanks{Department of Mathematics, Dartmouth College, Hanover, NH 03755, USA 
(\email{Lizuo.Liu@dartmouth.edu}, 
\email{Anne.E.Gelb@dartmouth.edu}
)}
\and
Lu Zhang\thanks{Department of Computational Applied Mathematics and Operations Research, and Ken Kennedy Institute, Rice University, Houston, TX 77005, USA}
(\email{lz82@rice.edu})
\and 
Anne Gelb\footnotemark[2]
}

\usepackage{amsopn}

\usepackage{comment}

\usepackage[normalem]{ulem}

\newcommand{\bmu}{{\boldsymbol\mu}}
\newcommand{\btheta}{{\boldsymbol\theta}}

\newcommand{\bu}{\mathbf{u}}

\newcommand{\cF}{\mathcal{F}}

\usepackage[normalem]{ulem}

\numberwithin{equation}{section}
\numberwithin{figure}{section}
\numberwithin{table}{section}

\allowdisplaybreaks

\patchcmd\newpage{\vfil}{}{}{}
\begin{document}

\maketitle

\begin{abstract}
	 We propose a neural entropy-stable conservative flux form neural network (NESCFN) for learning hyperbolic conservation laws and their associated entropy functions directly from solution trajectories, without requiring any predefined numerical discretization. While recent neural network architectures have successfully integrated classical numerical principles into learned models, most rely on prior knowledge of the governing equations or assume a fixed discretization. Our approach removes this dependency by embedding entropy-stable design principles into the learning process itself, enabling the discovery of physically consistent dynamics in a fully data-driven setting. By jointly learning both the numerical flux function and a corresponding entropy, the proposed method ensures conservation and entropy dissipation, critical for long-term stability and fidelity in the system of hyperbolic conservation laws. Numerical results demonstrate that the method achieves stability and conservation over extended time horizons and accurately captures shock propagation speeds, even without oracle access to future-time solution profiles in the training data.
\end{abstract}

\begin{keywords}
   Hyperbolic conservation laws, entropy stability, data-driven modeling, structure-preserving neural networks, conservative flux form neural networks
\end{keywords}

\begin{AMS}
        65M08, 68T07, 65M22, 65M32, 65D25
\end{AMS}



\section{Introduction}
\label{sec:introduction} 
Hyperbolic partial differential equations (PDEs) play a central role in modeling the dynamics of physical and engineered systems in which wave propagation and transport phenomena dominate. These equations govern the evolution of conserved quantities, such as mass, momentum, and energy, and are foundational to a broad range of applications in geophysical fluid dynamics, including ocean circulation, atmospheric flows, and sea ice mechanics (see, e.g.,~\cite{iceberg, Girard11, Euler, LeVeque02, Elastic-PlasticSeaIce, Shallow_water}). In many of these contexts, even when the precise form of the governing equations is unknown, it remains reasonable to assume that the underlying dynamics are governed by hyperbolic conservation laws. This observation motivates the development of data-driven approaches capable of inferring the governing structure of such systems directly from observed behavior, thereby enabling the prediction of future states based on past data.

With the growing availability of high-resolution observational and simulation data, machine learning has emerged as a powerful tool for modeling dynamical systems from data. Initial advances focused on learning ordinary differential equations (ODEs) from time series data (\cite{dynonet,rnn_ode,cdeeponet}), and recent efforts have extended these ideas to spatio-temporal systems governed by PDEs (\cite{neuralPDE,deeppropnet,yin2023continuous}). In the context of PDE learning, many existing methods are inspired by the classical method-of-lines framework. These approaches use neural networks to approximate the temporal derivative in a semi-discrete system, treating it as a black-box function of the state of the system at previous time steps; time integration is then performed using standard numerical solvers (see, e.g., \cite{chenPDE,Churchill_2023}). Often categorized as \emph{purely data-driven}, these methods do not incorporate explicit spatial structure or physical constraints. While they may succeed in short-term predictions, they often struggle to generalize or remain stable over long time horizons, especially in the presence of shocks, discontinuities, or conserved quantities.

To improve the fidelity of data-driven PDE models, several approaches have sought to embed spatial structure directly into the learning process. This is commonly done by incorporating spatial operators such as gradients, divergences, or Laplacians into the network architecture (e.g., \cite{kast2024positional,long2018pdenetlearningpdesdata}). These spatially-informed models offer improved stability and better representation of local interactions and geometric structure. However, despite these advances, such models typically do not enforce fundamental physical laws such as conservation or entropy dissipation, which are critical to the accuracy and reliability of solutions in many physical systems. The absence of these constraints can lead to non-physical artifacts and degraded performance over long-term integration.

To further bridge this gap, several methods have recently emerged that integrate principles from classical numerical schemes for conservation laws directly into neural network architectures. These approaches are specifically designed to capture the long-term dynamics of hyperbolic systems in a physically consistent way; see \cite{chencfn,liu2024entropy,Dimitrios2025gorinn,Seung2025tvdnet,RoeNet}. For instance, Godunov-Riemann informed neural networks introduced in \cite{Dimitrios2025gorinn} learn physical flux functions using numerical analysis-informed shallow neural networks, explicitly designed to preserve conservation laws and capture key wave interactions through embedded Riemann-solver structure. In contrast, inspired by total variation diminishing methods, the framework in~\cite{Seung2025tvdnet} incorporates a neural network closure into the discretized equations to suppress spurious oscillations and enforce boundedness to recovers hyperbolic phenomena by solving a PDE constrained optimization problem. 

Another prominent example is the conservative flux form neural network (CFN) introduced in~\cite{chencfn}, which mimics the structure of finite volume schemes by training a neural network to approximate the numerical flux function. This enables the recovery of unknown hyperbolic conservation laws purely from trajectory data. Similarly, RoeNet (\cite{RoeNet}) builds upon the classical Roe scheme, which approximates solutions to Riemann problems (\cite{roe}), embedding templaterizable Roe modules as data-driven computational primitives to predict the evolution of hyperbolic systems. Both CFN and RoeNet are trained to extrapolate system states using trajectory data observed in an initial time interval. However, neither method provides formal guarantees of entropy stability, which is critical for accurately modeling hyperbolic conservation laws. While~\cite{RoeNet} claims that entropy is not violated, no formal analysis or empirical validation is provided. Addressing this limitation, \cite{liu2024entropy} proposed the entropy-stable CFN (ESCFN), which explicitly encodes entropy stability into the neural network architecture in the form of slope-limiting numerical methods. Specifically, their method employs the second-order accurate, non-oscillatory Kurganov–Tadmor (KT) scheme (\cite{KTscheme2000}) and trains a neural network to approximate the numerical flux term within this framework, thereby ensuring entropy stability in the learned model.

While the ESCFN framework proposed in~\cite{liu2024entropy} is both simple and empirically effective for predicting the dynamics of hyperbolic conservation laws, its design relies on a pre-specified numerical scheme (e.g., KT discretization) as a foundation. This reliance also implicitly assumes a specific class of corresponding entropy pairs (discussed in more detail in \Cref{sec: preliminaries}), thereby limiting the method's flexibility in scenarios where the governing equations are entirely unknown or difficult to approximate using classical methods. To overcome this limitation, we draw inspiration from the entropy-stable numerical flux design principles developed in~\cite{tadmor2016entropy} and propose a new data-driven framework that dispenses with the need for any predefined numerical scheme. Our approach simultaneously learns both the underlying hyperbolic conservation law and a corresponding entropy function directly from solution trajectories. This enables the identification of physically meaningful and dynamically consistent models purely from data. The resulting framework preserves essential structural properties, such as conservation and entropy dissipation, while offering enhanced generality and adaptability across a broad class of systems.

The rest of this paper is organized as follows. In \Cref{sec: preliminaries} we review the key theoretical and numerical foundations relevant to hyperbolic systems of conservation laws.
\Cref{sec: learned entropy} introduces neural entropy-stable conservative flux form neural networks (NESCFN) for learned hyperbolic conservation laws directly from solution trajectories. The experimental setup and evaluation metrics are described in \Cref{sec: training design}. In \Cref{sec: numerical}, we present a series of numerical examples that demonstrate the effectiveness and robustness of the proposed NESCFN framework. Finally, we conclude with a summary and briefly discuss directions for future research in \Cref{sec:conclusion}.

\section{Preliminaries}
\label{sec: preliminaries} We first review the essential theoretical and numerical foundations associated with hyperbolic systems of conservation laws. These include basic definitions, entropy conditions, the symmetrization framework, and the classical numerical schemes, which serve both to establish the notation and to contextualize the key principles underlying our later development.

\subsection{Hyperbolic conservation laws and entropy pairs}

 Our primary goal is to study the behavior of solutions to a hyperbolic system of conservation laws, especially in contexts where the governing equations are not known in closed form. To this end, we briefly review the structural properties of such systems, which are essential for the construction of reliable numerical schemes.
 
 Consider the general $d$-dimensional system of conservation laws written as 
\begin{equation} \label{eq: conservation_law} \frac{\partial \mathbf{u}}{\partial t}  + \sum_{j=1}^d \frac{\partial \mathbf{f}_j(\mathbf{u})}{\partial x_j}  = \mathbf{0}, \quad \mathbf{x} = [x_1, \cdots, x_d]^\top \in \Omega \subset \mathbb{R}^d, \quad t \in (0, T), 
\end{equation} 
subject to appropriate initial and boundary conditions. Here, $\mathbf{u} = [u^1, \ldots, u^p]^\top$ represents the vector of conserved state variables, taking values in a convex set $\mathcal{D} \subset \mathbb{R}^p$. The flux functions $\mathbf{f}_j : \mathcal{D} \rightarrow \mathbb{R}^p$ are assumed to be sufficiently smooth.

For each spatial direction $1 \leq j \leq d$, define the flux Jacobian matrix 
\begin{equation}\label{eq: flux_jacobi}
A_j({\bf u}) := \mathbf{f}_j'(\mathbf{u})=\left\{\frac{\partial f_j^i(\mathbf{u})}{\partial u^k}\right\}_{1 \leq i, k \leq p}.
\end{equation}
The system \eqref{eq: conservation_law} is called \emph{hyperbolic} if for any unit vector $\boldsymbol{\nu} = [\nu_1,\cdots, \nu_d]^\top \in \mathbb{R}^d$, the matrix $A(\mathbf{u}, \boldsymbol{\nu}) := \sum_{j=1}^d \nu_j A_j(\mathbf{u})$ has $p$ real eigenvalues and a complete set of linearly independent eigenvectors for all $\mathbf{u} \in \mathcal{D}$.

It is well known that even when the initial data are smooth, solutions to hyperbolic systems can develop discontinuities, such as shocks or contact discontinuities, within finite time. As such, solutions must be interpreted in the weak sense. Weak solutions, however, are not necessarily unique. To enforce physical admissibility and uniqueness, one typically imposes additional constraints in the form of entropy conditions.

\subsubsection{Entropy pairs and admissibility} Entropy functions are introduced to identify admissible weak solutions. These are convex scalar functions $\eta({\bf u}): \mathcal{D} \rightarrow \mathbb{R}$, paired with corresponding entropy fluxes ${G}_j({\bf u}): \mathcal{D} \rightarrow \mathbb{R}, 1 \leq j \leq d$, satisfying the compatibility condition given by \cite{godlewski2013numerical}: 
\begin{definition}
Suppose that the domain $\mathcal{D}$ is convex. A convex function $\eta(\mathbf{u})$ is called \emph{entropy function} for the system \eqref{eq: conservation_law} if there exist associated entropy fluxes $G_j(\mathbf{u})$, $1 \leq j \leq d$, such that
\[\eta'(\mathbf{u}) \mathbf{f}_j'(\mathbf{u})=G_j^{\prime}(\mathbf{u})\quad \forall j,\]
where the derivatives are understood as row vectors.
\end{definition}

In regions where the solution is smooth, multiplying \eqref{eq: conservation_law} from the left by $\eta'(\mathbf{u})$ provides an additional conservation law
\[\frac{\partial \eta(\mathbf{u})}{\partial t}+\sum_{j=1}^d\frac{\partial G_j(\mathbf{u})}{\partial x_j}=0.\]
At discontinuities, however, the entropy must dissipate. This leads to the following entropy inequality condition \eqref{eq:entropy_cond}, which characterizes physically admissible solutions:
\begin{definition}
 A weak solution $\mathbf{u}$ of \eqref{eq: conservation_law} is called an entropy solution if for any entropy function $\eta({\bf u})$ the following inequality holds in the sense of distribution: \begin{equation} \label{eq:entropy_cond} 
\frac{\partial \eta(\mathbf{u})}{\partial t} + \sum_{j=1}^d \frac{\partial G_j(\mathbf{u})}{\partial x_j} \leq 0. \end{equation}
\end{definition}

\subsubsection{Symmetrization and entropy variables} 
The central task is then to identify entropy functions for systems of conservation laws \eqref{eq: conservation_law}. In the scalar case ($p=1$), any convex function qualifies as an entropy. However, for systems with $p>1$, identifying entropy functions becomes much more difficult. The following result, adapted from \cite{godlewski2013numerical}, characterizes the existence of entropy functions by the symmetry of certain coefficient matrices.

\begin{theorem}\label{thm: existence entropy}
Let $\eta: \mathcal{D} \rightarrow \mathbb{R}$ be a strictly convex function. Then $\eta$ is an entropy for the system \eqref{eq: conservation_law} if and only if the matrices $\eta''({\bf u}){\bf f}_j'({\bf u}), 1\leq j\leq d$, are symmetric.
\end{theorem}

Moreover, the existence of a strictly convex entropy function $\eta$ ensures that the system \eqref{eq: conservation_law} is symmetrizable and thus hyperbolic. In fact, since $\eta$ is strictly convex, a one-to-one change of variables can be defined by the \emph{entropy variable} ${\bf v} = (\eta'({\bf u}))^\top$. Letting 
\[{\bf g}_j({\bf v}) := {\bf f}_j({\bf u}({\bf v})),\]
the conservation law \eqref{eq: conservation_law} can be reformulated as
\begin{equation}\label{eq: symmetric}
\mathbf{u}'(\mathbf{v}) \frac{\partial \mathbf{v}}{\partial t}+\sum_{j=1}^d \mathbf{g}_j'(\mathbf{v}) \frac{\partial \mathbf{v}}{\partial x_j}=0.
\end{equation}
Given that $\eta({\bf u})$ is strictly convex and $\eta''({\bf u}){\bf f}_j'({\bf u})$ are symmetric, we compute
\begin{equation*}
\begin{aligned}
{\bf g}_j'({\bf v}) &= {\bf f}_j'({\bf u}){\bf u}'({\bf v}) = {\bf f}_j'({\bf u})(\eta''({\bf u}))^{-1} \\
&= (\eta''({\bf u}))^{-1}({\bf f}_j'({\bf u}))^\top \eta''({\bf u})(\eta''({\bf u}))^{-1} = (\eta''({\bf u}))^{-1}({\bf f}_j'({\bf u}))^\top,
\end{aligned}
\end{equation*}
where the second equality follows from the definition of the entropy variable, and the third equality results from the symmetry of $\eta''({\bf u})$. Hence, we have ${\bf g}_j'({\bf v}) = ({\bf g}_j'({\bf v}))^\top$, indicating that each ${\bf g}_j'({\bf v})$ is symmetric. Therefore, the system \eqref{eq: symmetric} is in \emph{symmetrized form}. Furthermore, the flux Jacobian of the system \eqref{eq: conservation_law} in an arbitrary direction $\boldsymbol{\nu}\in \mathbb{R}^d$ is given by $A(\mathbf{u}, \boldsymbol{\nu})=\sum_{j=1}^d \nu_j \mathbf{f}_j'(\mathbf{u})=\sum_{j=1}^d \nu_j \mathbf{g}_j'(\mathbf{v}) \mathbf{v}'(\mathbf{u})$, which is similar to the symmetric matrix $\mathbf{v}'(\mathbf{u})^{\frac{1}{2}}\big(\sum_{j=1}^d \nu_j \mathbf{g}_j'(\mathbf{v})\big) \mathbf{v}'(\mathbf{u})^{\frac{1}{2}}.$
Hence, the existence of a strictly convex entropy function implies the hyperbolicity of the system \eqref{eq: conservation_law}.

This symmetrization result plays a fundamental role in both the theoretical and computational treatment of hyperbolic conservation laws. In particular, it establishes a direct link between entropy structure and hyperbolicity, which is crucial for the design of well-posed and stable numerical schemes. With this theoretical foundation in place, we now turn our attention to a review of numerical schemes for hyperbolic conservation laws. This review includes both classical approaches and recent data-driven advances, and serves as a basis for understanding and motivating the NESCFN introduced in Section \ref{sec: learned entropy}.




\subsection{Numerical schemes for conservation laws}\label{sec: review numerical}
For illustrative purposes, we restrict our attention to the one-dimensional case, corresponding to $d = 1$ in \eqref{eq: conservation_law}. This simplification is made for ease of presentation, but we emphasize that our proposed framework naturally extends to multidimensional systems. Indeed, we provide numerical results both in 1D and 2D in Section \ref{sec: numerical} to demonstrate the broader applicability of the proposed method.

The one-dimensional conservation law under consideration is given by
\begin{equation}\label{eq:conservation_law_1d}
\frac{\partial \mathbf{u}}{\partial t} + \frac{\partial \mathbf{f}(\mathbf{u})}{\partial x} = 0, \quad x \in \Omega = (a,b), \quad t \in (0,T), 
\end{equation} 
where ${\bf u}\in \mathbb{R}^p$ represents the conserved states, and ${\bf f}({\bf u}) \in \mathbb{R}^p$ is the corresponding flux function.

To numerically solve the problem \eqref{eq:conservation_law_1d}, we partition the computational domain $\Omega$ into uniform spatial grid points $\{x_j\}_{j=0}^n$, where $x_j= j\Delta x, \Delta x = \frac{b-a}{n}.$
We note that the notation $x_j$ is reused here to denote the $j$-th spatial grid point, which differs from its earlier use in \eqref{eq: conservation_law}, where it referred to the $j$-th spatial direction in multidimensional problems. This slight abuse of notation is standard in numerical schemes and should not cause confusion in the 1D setting considered here. We further approximate the solution by cell averages ${\bf u}_j(t) = [u_j^1(t), \cdots, u_j^p(t)]^T$ over the cell $I_j := (x_j - \frac{\Delta x}{2}, x_j + \frac{\Delta x}{2})$ with 
\begin{equation}\label{eq: cell_average}
u_j^i(t)=\int_{I_j} u^i(x, t) d x, \quad j=1, \cdots, n-1; \ i = 1,\cdots,p,
\end{equation}
and suitable boundary conditions are enforced at the endpoints denoted by $u_0^i$ and $u_n^i$. Note that, also with a slight abuse of notation, we use ${\bf u} = [u^1, \cdots, u^p]^\top$ to denote the continuous solution and ${\bf u}_j = [u_j^1, \cdots, u_j^p]^\top$ to denote its spatially discretized counterpart for the cell average over $I_j$. This convention will be used throughout the paper unless otherwise specified, and any potential ambiguity will be explicitly clarified when necessary.

Within this framework, we now review two critical classes of schemes: entropy-conservative and entropy-stable schemes, developed in \cite{Tadmor1987EntropyStable, tadmor2016entropy}.

\subsubsection{Entropy-conservative schemes} We start with a class of second-order entropy-conservative schemes introduced in \cite{Tadmor1987EntropyStable}. Let $v_j^i(t), 1\leq i\leq p$, be the approximation to the entropy variable $v^i(x,t)$ averaged over the cell $I_j$ according to the rule in \eqref{eq: cell_average}, and denote ${\bf v}_j = [v_j^1, \cdots, v_j^p]^\top$. The second-order entropy conservative scheme is then given by:
\begin{theorem}\label{Thm:conservative}
Let $(\eta({\bf v}), G({\bf v}))$ be an entropy pair associated with the conservation law \eqref{eq:conservation_law_1d}. Then, the numerical scheme
\begin{equation}\label{eq:numerical_scheme}
\frac{d}{dt} \left[ {\bf u}({\bf v}_j(t)) \right] = -\frac{1}{\Delta x} \left[ {\bf g}^*_{2,j+1/2} - {\bf g}^*_{2,j-1/2} \right]
\end{equation}
is entropy conservative, with the numerical flux 
\begin{equation}\label{eq:flux_entropy_conservative}
{\bf g}^*_{2,j+1/2}({\bf v}_j, {\bf v}_{j+1}) = \int_{0}^{1} {\bf g}\big({\bf v}_j + \xi ({\bf v}_{j +1} - {\bf v}_j)\big) d\xi,
\end{equation}
where the integral of the vector function ${\bf g}$ is understood as component-wise. In addition, the scheme \eqref{eq:numerical_scheme} satisfies the following cell entropy equality
\begin{equation}
\frac{d}{dt} \eta({\bf v}_j(t)) + \frac{1}{\Delta x} \left( G^*_{2,j+1/2} - G^*_{2,j-1/2} \right) = 0
\label{eq: entropy conservative scheme}
\end{equation}
with the numerical entropy flux is given by
\begin{equation}\label{eq:G}
\begin{aligned}
    G^*_{2,j+1/2} & = \frac{1}{2} ({\bf v}_j + {\bf v}_{j+1})^\top {\bf g}^*_{2,j+1/2} + \frac{1}{2} \big(G({\bf v}_j) + G({\bf v}_{j+1})\big)\\
    &- \frac{1}{2} \big({\bf v}_j^\top {\bf g}({\bf v}_j) + {\bf v}_{j+1}^\top {\bf g}({\bf v}_{j+1})\big). 
    \end{aligned}
\end{equation}

\end{theorem}

We note that the second-order entropy-conservative numerical flux \eqref{eq:flux_entropy_conservative} can be extended to construct higher order schemes. Following the formulation in \cite{LeFlochEntropyConservativeSchemeArbitrayOrder}, one can define entropy-conservative numerical fluxes of arbitrary even order $2m, m\in\mathbb{N}^+$, by appropriate linear combinations:
\begin{equation}\label{eq: conservative_flux_high_order}
\begin{aligned}
& {\bf g}_{2m, j+1/2}^\ast({\bf v}_{j-m+1}, \cdots, {\bf v}_{j+m})\\
= &\sum_{i = 1}^m \alpha_{i,m} \left({\bf g}^\ast_{2,j+1/2}({\bf v}_j,{\bf v}_{j+i}) + \cdots + {\bf g}_{2,j+1/2}^\ast({\bf v}_{j-i+1},{\bf v}_{j+1})\right),
\end{aligned}
\end{equation}
with coefficients $\alpha_{i,m}$ satisfying the following moment conditions
\[2\sum_{i=1}^m i\alpha_{i,m} = 1, \quad \sum_{i=1}^m i^{2s-1}\alpha_{i,m} = 0 \ \text{ for } s = 2,\cdots, m.\]

\begin{remark}
Given the relation ${\bf g}({\bf v}) = {\bf f}\big({\bf u}({\bf v})\big)$, the entropy-conservative numerical scheme \eqref{eq:numerical_scheme} can be equivalently written as
\begin{equation}\label{eq: numerical_scheme_u}
\frac{d}{dt} {\bf u}_j(t) = -\frac{1}{\Delta x} \left[{\bf f}_{j+1/2}^\ast - {\bf f}_{j-1/2}^\ast\right] \quad \text{with}\quad {\bf f}_{j+1/2}^\ast = {\bf g}_{j+1/2}^\ast
\end{equation}
in terms of conservative states ${\bf u}$ and its corresponding flux ${\bf f}({\bf u})$. In addition, for notational simplicity, we omit the order subscript $2m$ in the numerical fluxes ${\bf f}_{2m,j+1/2}^\ast$ and ${\bf g}_{2m,j+1/2}^\ast$ from this point onward. 

\end{remark}

\subsubsection{Entropy-stable schemes}\label{sec: review_entropy_stable}
While entropy-conservative schemes preserve the entropy of the conservation law without introducing numerical dissipation, they do not necessarily suppress spurious oscillations, especially in the presence of shocks or discontinuities. In practical applications, ensuring entropy stability is essential to maintain the robustness and physical admissibility of the numerical solution. A straightforward approach to constructing entropy-stable schemes by augmenting entropy-conservative schemes with appropriate dissipation terms was proposed in \cite{tadmor2016entropy}, and is outlined in the theorem below.
\begin{theorem}\label{TH:entropy_stable}
   Let $(\eta({\bf v}), G({\bf v}))$ be an entropy pair associated with the conservation law \eqref{eq:conservation_law_1d}, and let ${\bf f}_{j+1/2}^\ast$ denote an entropy-conservative flux as in \eqref{eq: numerical_scheme_u} and $G^\ast_{j+1/2}$ be the numerical entropy flux defined in \eqref{eq:G}. Suppose that $D_{j+1/2} \succeq 0$ is a symmetric positive semidefinite matrix. Then, the scheme 
    \begin{equation}\label{eq: numerical_scheme_s}
    \frac{d}{dt} {\bf u}_j(t) = -\frac{1}{\Delta x} \left(\hat {\bf f}_{j+1/2} - \hat {\bf f}_{j-1/2}\right)
    \end{equation}
with the numerical flux
\begin{equation}
\hat{\bf f}_{j+1/2} = {\bf f}_{j+1/2}^\ast - \frac{1}{2} D_{j+1/2} [[{\bf v}]]_{j+1/2}
\label{eq: entropy stable flux}
\end{equation}
is entropy stable, satisfying
\[\frac{d}{dt} \eta({\bf v}_j(t)) + \frac{1}{\Delta x} \left( \hat G_{j+1/2} - \hat G_{j-1/2} \right) \leq 0.\]
Here, $[[{\bf v}]]_{j+1/2} := {\bf v}_{j+1} - {\bf v}_j $ denotes the jump of the entropy $\bf v$ at the cell edge $x_{j+1/2}$, and the numerical entropy flux function $\hat G_{j+1/2}$ takes the following form
\[\hat G_{j+1/2} = G^\ast_{j+1/2} + \frac{1}{2}\bar {\bf v}_{j+1/2} D_{j+1/2} [[{\bf v}]]_{j+1/2} \quad \text{with}\quad \bar{\bf v}_{j+1/2} = \frac{1}{2}({\bf v}_j + {\bf v}_{j+1}).\]
\end{theorem}

We note that the construction of entropy-stable schemes depends on the appropriate choice of the diffusion matrix $D_{j+1/2}$. While many choices exist, in this paper we focus on a specific form related to the flux Jacobian of the underlying hyperbolic conservation law. Details are given in Section \ref{sec: learned entropy}.

It should also be emphasized that the numerical flux defined in \eqref{eq: entropy stable flux} generally gives only spatial first-order accuracy, a limitation arising from the fact that the jump term $[[\mathbf{v}]]_{j+1/2}$ is intrinsically $\mathcal{O}(\Delta x)$. Notably, this remains true even if the entropy-conservative flux in \eqref{eq:flux_entropy_conservative} is replaced by a higher order version, such as the high-order entropy-conservative flux given in \eqref{eq: conservative_flux_high_order}. Achieving higher order accuracy requires more than simply upgrading the flux function, it requires a refined reconstruction of the entropy variables $\mathbf{v}$ within each cell $I_j$. Specifically, one must use a higher degree polynomial reconstruction rather than the standard piecewise constant approximation in order for the jump $[[\mathbf{v}]]_{j+1/2}$ to achieve high-order accurate. For example, \eqref{eq: upm} provides a piecewise linear polynomial reconstruction to obtain a second-order method. However, due to space limitations, we omit the detailed formulation of such high-order entropy-stable fluxes, and refer interested readers to the comprehensive treatment in \cite{fjordholm2012arbitrarily}, which discusses the design and analysis of arbitrarily high-order accurate entropy-stable schemes. 

\subsubsection{Data-driven extensions}\label{sec:KT_CFN}

Recent advances in scientific machine learning have opened new avenues for integrating data-driven modeling with classical numerical methods. One promising direction is the incorporation of conservation law structures, such as entropy conservation and stability, into neural network-based solvers.  A representative example is the framework proposed in \cite{liu2024entropy}, which introduced a family of entropy-stable conservative flux-form neural networks (CFNs). These models embed entropy stability directly into the architecture of the neural network by extending classical finite-volume schemes, specifically, the second-order, non-oscillatory Kurganov--Tadmor (KT) method (\cite{kurganov2000new}). Unlike traditional schemes that require explicit knowledge of the flux function, the CFN learns the flux from the data while preserving essential mathematical structures such as conservation and entropy dissipation. The approach is briefly reviewed below.

\subsubsection*{Kurganov--Tadmor entropy stable CFN} The KT-ESCFN retains the conservative form of the original KT scheme, but replaces the known flux with a trainable neural network. The semi-discrete update rule is written as 
\begin{equation*}
\frac{d}{d t} {\bf u}_j(t)=-\frac{H^{NN}_{j+1/2}(t)-H^{NN}_{j-1 / 2}(t)}{\Delta x},
\label{eq: Kurganov Tadmor Scheme}
\end{equation*}
where $H^{NN}_{j+\frac{1}{2}}$ is the learned numerical flux defined by
\begin{equation*}
H^{NN}_{j+1 / 2}(t)=\frac{\cF^{\bmu}\left({\bf u}_{j+1 / 2}^{+}(t)\right)+\cF^{\bmu}\left({\bf u}_{j+1 / 2}^{-}(t)\right)}{2}-\frac{a^{NN}_{j+1 / 2}(t)}{2}\left[{\bf u}_{j+1 / 2}^{+}(t)-{\bf u}_{j+1 / 2}^{-}(t)\right].
\label{eq: Numerical flux NN}
\end{equation*}
Here, $\mathcal{F}^{\bmu}$ denotes a fully connected feed-forward neural network parameterized by trainable weights $\bmu$, which replaces the analytical flux function ${\bf f}({\bf u})$ in the original KT scheme (see \cite{kurganov2000new}). The local wave propagation speed $a^{NN}_{j+1/2}(t)$ is also learned from the data using a second neural network $\rho_{\bf w}$ that approximates the spectral radius of the Jacobian with 
\begin{equation}
\alpha^{NN}_{j+1 / 2} = \max \left\{\left|\rho_{\bf w}\left(\frac{\partial \cF^{\bmu}}{\partial {\bf u}}\left({\bf u}_{j+1 / 2}^{+}(t)\right)\right)\right|, \left|\rho_{\bf w}\left(\frac{\partial \cF^{\bmu}}{\partial {\bf u}}\left({\bf u}_{j+1 / 2}^{-}(t)\right)\right)\right|\right\},
\label{eq: rho_W}
\end{equation}
 where $\rho_{\bf w}$ is a neural network mapping the full Jacobian matrix to a positive scalar approximation of its dominant eigenvalue. Additionally, the states ${\bf u}_{j+1/2}^{\pm}$ are defined as in the original KT scheme (\cite{kurganov2000new}) with
\begin{equation}
{\bf u}_{j+1 / 2}^{+}(t)={\bf u}_{j+1}(t)-\frac{\Delta x}{2}\left({\bf u}_x\right)_{j+1}(t), \quad {\bf u}_{j+1 / 2}^{-}(t)={\bf u}_j(t)+\frac{\Delta x}{2}\left({\bf u}_x\right)_j(t),
\label{eq: upm}
\end{equation}
where the spatial derivatives $({\bf u}_x)_j$ are computed using a Total Variation Diminishing (TVD) limiter. In particular, for each component $(u_x^i)_j, 1\leq i \leq p$, of $({\bf u}_x)_j$ we have
\begin{equation*}
( {u}_{x}^i )_{j} = \psi({r}) \left( {u}_{j + 1}^i -  {u}_{j}^i \right), \quad {r} = \dfrac{{u}_{j}^i - {u}_{j-1}^i}{{u}_{j+1}^i - {u}_{j}^i},
\end{equation*}
with the minmod-type limiter $\psi(r) = \max\big( 0, \min( r, (1+r)/2, 1) \big).$

This formulation ensures that the learned flux retains the entropy-stable structure of the KT scheme, while allowing data-driven discovery of flux functions and wave speed. The KT-ESCFN achieves strong performance in approximating the solution of hyperbolic systems, exhibiting stability and accuracy even under noisy and sparsely sampled observations. Notably, these models remain robust even when trained on data that do not contain discontinuities, demonstrating their ability to generalize to complex regimes governed by underlying physical laws.

\subsubsection*{Motivation for a new approach}  While the KT-ESCFN framework of \cite{liu2024entropy} marks a significant step forward in combining neural networks with structure-preserving numerical schemes, it still requires the prior selection of a classical discretization framework (e.g.,~the KT scheme) as scaffolding for the learning process. In other words, the numerical structure must be manually specified in advance, which limits the flexibility of the model in scenarios where the most appropriate discretization is unknown or system-dependent. 

In this work, we build on insights from \cite{liu2024entropy} and the theory of entropy-stable numerical methods to propose a new learning-based framework, NESCFN. Our approach seeks to relax the dependence on predefined schemes while preserving the essential entropy properties through learnable components. The details of this new methodology are presented next.

\section{Entropy-stable network with neural entropy}\label{sec: learned entropy} 
As discussed in Section \ref{sec: review_entropy_stable}, the construction of entropy-stable numerical schemes for hyperbolic conservation laws relies on augmenting entropy-conservative fluxes with appropriate dissipation. According to Theorem \ref{TH:entropy_stable}, such schemes can be formulated by adding a numerical viscosity term to an entropy-conservative flux, and can be constructed using symmetric positive semidefinite matrices together with the jump of the entropy variables across adjacent cells, see \cref{eq: entropy stable flux}. This general strategy provides a flexible basis for ensuring entropy stability.

While Theorem~\ref{TH:entropy_stable} allows any symmetric positive semidefinite matrix $D_{j+1/2}$ in the formulation, practical implementations often make specific choices that are consistent with the structure of the underlying conservation law. A common approach is to derive the diffusion matrix based on the eigendecomposition of the flux Jacobian (see, e.g., \cite{fjordholm2012arbitrarily,FarzadAffordableEntropyConsistentEulerFlux}): Let ${\bf u}_{j+1/2}$ denote the numerical approximation of the state variable ${\bf u}$ evaluated at the cell edge $x_{j+1/2}$, and $A_{j+1/2} = {\bf f}'({\bf u}_{j+1/2})$ be the flux Jacobian matrix associated with the conservation law \eqref{eq:conservation_law_1d} evaluated at ${\bf u}_{j+1/2}$. Due to the hyperbolicity of the system, $A_{j+1/2}$ has a complete set of linearly independent eigenvectors and can be diagonalized as
\[A_{j+1/2} = R_{j+1/2} \Lambda_{j+1/2} R_{j+1/2}^{-1},\]
where $\Lambda_{j+1/2} = \text{diag}(\lambda_{j+1/2}^1, \ldots, \lambda_{j+1/2}^p)$ is the diagonal matrix of the eigenvalues, and $R_{j+1/2}$ is the matrix of the corresponding right eigenvectors. A commonly used entropy-stable diffusion operator is the Rusanov matrix, defined as
\begin{equation}
D_{j+1/2} = \tilde{R}_{j+1/2} (\lambda_{j+1/2}^{\max} \mathbb{I}_p) \tilde{R}^\top_{j+1/2},
\label{eq: entropy D mat}
\end{equation}
where $\lambda_{j+1/2}^{\max} = \max(|\lambda_{j+1/2}^1|, \ldots, |\lambda_{j+1/2}^p|)$ is the maximum wave speed, $\mathbb{I}_p \in \mathbb{R}^{p \times p}$ is the $p\times p$ identity matrix, and $\tilde{R}_{j+1/2}$ is a rescaled version of $R_{j+1/2}$ such that the entropy Hessian satisfies $\mathbf{u}'({\bf v}_{j+1/2}) = \tilde{R}_{j+1/2} \tilde{R}_{j+1/2}^\top$. For a rigorous discussion of the existence of such $\tilde{R}_{j+1/2}$, we refer the reader to Theorem 4 of \cite{Barth1999NumericalMethodsGasDynamic}.

Substituting $D_{j+1/2}$ in \eqref{eq: entropy D mat} into the general entropy-stable flux formula from Theorem \ref{TH:entropy_stable} yields the following Rusanov-type entropy-stable numerical flux:
\begin{equation}
\begin{aligned}
    \hat{\bf f}_{j+1/2} &= {\bf f}^\ast_{j+1/2} - \frac{1}{2} \tilde{R}_{j+1/2} \big(\lambda_{j+1/2}^{\max}\mathbb{I}_p\big) \tilde{R}^\top_{j+1/2}  [[{\bf v}]]_{j+1/2}\\
    &= {\bf f}^\ast_{j+1/2} - \frac{1}{2} \lambda_{j+1/2}^{\max}\tilde{R}_{j+1/2}  \tilde{R}^\top_{j+1/2}  [[{\bf v}]]_{j+1/2} \\
    &= {\bf f}^\ast_{j+1/2} - \frac{1}{2}\lambda_{j+1/2}^{\max}{\bf u}'({\bf v}_{j+1/2})[[{\bf v}]]_{j+1/2} \\
    &= {\bf f}^\ast_{j+1/2} - \frac{1}{2} \lambda_{j+1/2}^{\max}\big(\eta''({\bf u}_{j+1/2})\big)^{-1} [[\eta'({\bf u})]]_{j+1/2},
\end{aligned}
\label{eq: Rusanov entropy stable flux}
\end{equation}
where $\mathbf{f}^\ast_{j+1/2}$ is the entropy-conservative numerical flux (defined earlier in \eqref{eq: numerical_scheme_u}), and $[[\cdot ]]_{j+1/2}$ denotes the jump of the quantity across the cell edge $x_{j+1/2}$.

This Rusanov-type flux is widely adopted due to its simplicity, inherent entropy stability, and compatibility with a variety of system structures, see \cite{fjordholm2012arbitrarily,FarzadAffordableEntropyConsistentEulerFlux}. In what follows, we build upon this formulation to design an entropy-stable conservative flux neural network with neural entropy for learning hyperbolic conservation laws.

\subsection{Learning entropy-stable flux via neural parameterization}\label{sec:learning entropy stable}

To incorporate the entropy structure directly into a learnable framework, we first parameterize the entropy function using an input convex neural network (see Section \ref{sec:input_cnn} for details), denoted by $\eta_{\btheta}$ with $\btheta$ denote the trainable variables. This architecture guarantees convexity of the learned entropy by design, which is one of the critical properties that underpins the hyperbolicity of the resulting learned conservation law system. Once $\eta_{\btheta}$ is defined, both the parameterized entropy variable ${\bf v}_{\btheta} = \eta_{\btheta}'({\bf u})$ and the parameterized entropy Hessian $\eta_{\btheta}''({\bf u})$ can be easily computed via auto-differentiation.

 Given this parameterization, the entropy-stable flux can be formulated in a Rusanov-like structure, with the dissipation term expressed in terms of the learned entropy as\vspace{-0.25cm}
\begin{equation}
    \hat{\bf f}_{j+1/2} =  {\bf f}^\ast_{j+1/2} - \frac{1}{2} \lambda_{j+1/2}^{\max} \left(\eta_{\btheta}''({\bf u}_{j+1/2})\right)^{-1} [[\eta_\btheta'({\bf u})]]_{j+1/2}.
\label{eq: Param entropy conservative flux}
\end{equation}

To complete this formulation, both $\lambda^{\max}_{j+1/2}$ and the entropy-conservative flux $\mathbf{f}^\ast_{j+1/2}$ must also be learned. Inspired by the framework introduced in \cite{liu2024entropy} (see also Section \ref{sec:KT_CFN}), we approximate the entropy-conservative flux $\mathbf{f}^\ast_{j+1/2}$ by an average of neural flux predictions with
\begin{equation}\label{eq:entropy_stable_flux_nn}
{\bf f}_{j+1/2}^\ast \approx \mathcal{F}^{\bmu,\ast}_{j+1/2} = \frac{1}{2}\Big(\mathcal{F}^\bmu({\bf u}_{j+1/2}^+) + \mathcal{F}^\bmu({\bf u}_{j+1/2}^-)\Big),
\end{equation}
where $\mathcal{F}^\bmu$ is a fully connected feedforward neural network with trainable parameters $\bmu$, and $\mathbf{u}_{j+1/2}^\pm$ are reconstructed interface values (e.g., as in \eqref{eq: upm}). The local wave speed $\lambda^{\max}_{j+1/2}$ for the conservation law system is approximated by a neural surrogate, $\rho_{\bf w}( {\bf u}_{j+1/2})$\footnote{For scalar hyperbolic conservation laws, we use the absolute value of the Jacobian to approximate the local wave speed.}, designed to explicitly enforce the CFL condition with
 \[\rho_{\bf w}( {\bf u}_{j+1/2}) = \dfrac{\Delta x}{\Delta t}\Big[ 1 - \tanh\Big(\big| \phi_{\bf w}\big( (\mathcal{F}^\bmu)'({\bf u}_{j+1/2})\big) \big| \Big)\Big],\] 
 where $\phi_{\bf w}(\cdot)$ is a neural network parameterized by $\bf w$ that takes the Jacobian of $\mathcal{F}^\bmu$ as input and returns a positive scalar output. This design ensures bounded and positive wave speed approximations while maintaining numerical stability. The details of the structure of neural networks and the choice of training parameters are deferred to Section~\ref{sec:networkdetails}. 

Combining these elements, we define the final implementable neural entropy-stable flux as
\begin{equation}
\begin{aligned}
   \hat{\bf f}_{j+1/2} \approx \hat{\mathcal F}^{\bmu,{\bf w, \btheta}}_{j+1/2} = \mathcal{F}^{\bmu, \ast}_{j+1/2} - &\frac{1}{2}  \underbrace{\text{max}\left(\rho_{\bf w}\left( {\bf u}_{j+1/2}^{+} \right), \rho_{\bf w}\left( {\bf u}_{j-1/2}^{-} \right)\right)}_{ \approx  \lambda_{j+1/2}^{\max}} \\
    &\cdot \big(\eta_{\btheta}''(\bar {\bf u}_{j+1/2})\big)^{-1}\left(\eta'_\btheta({\bf u}_{j+1/2}^+) - \eta'_\btheta({\bf u}_{j+1/2}^-)\right),
    \end{aligned}
\label{eq: neural-entropy stable flux}
\end{equation}
where $\bar{\bf u}_{j+1/2} = \frac{1}{2}({\bf u}_{j+1/2}^- + {\bf u}_{j+1/2}^+)$ is the cell-averaged state at the interface. Finally, to ensure numerical stability when inverting the entropy Hessian in \eqref{eq: neural-entropy stable flux}, we apply Cholesky decomposition to a regularized form of the Hessian. Specifically, we add an epoch-dependent diagonal perturbation with exponentially decaying magnitude
\begin{equation}
 L({\bf u}_{j+1/2})L({\bf u}_{j+1/2})^\top = \eta_{\btheta}''({\bf u}_{j+1/2}) + C_1 C_2^{\text{Epoch}}\mathbb{I}_p, 
\label{eq: stable cho decomp}
\end{equation}
where \(C_1 > 0\) and \(0< C_2 < 1\) are fixed constants (see Section~\ref{sec:networkdetails}). The resulting linear system is then solved efficiently using routines such as \texttt{cho\_solve} in JAX to compute the inverse Hessian-vector product required in \eqref{eq: Param entropy conservative flux}.

Thanks to the structural guarantees of convex input networks and Theorem~\ref{TH:entropy_stable}, the numerical flux defined in \eqref{eq: neural-entropy stable flux} is entropy-stable for the learned conservation law 
\begin{equation} \label{eq: learned_conservation_law} 
\frac{\partial \mathbf{u}}{\partial t} + \frac{\partial \mathcal{F}^\bmu(\mathbf{u})}{\partial x} = 0
\end{equation} 
with respect to the learned entropy $\eta_\btheta$. Moreover, if the training procedure for the neural networks is designed to enforce the required symmetric structure in Theorem~\ref{thm: existence entropy}, then the learned conservation law is guaranteed to be hyperbolic. With the learned entropy-stable flux \eqref{eq: neural-entropy stable flux}, we have the following semidiscrete entropy-stable scheme for the learned conservation law \eqref{eq: learned_conservation_law}
\begin{equation}\label{eq: scheme_learned_law}
\frac{d}{dt} {\bf u}_j(t) = -\frac{1}{\Delta x} \left( \hat{\mathcal{F}}_{j+1/2}^{\bmu,{\bf w}, \btheta} - \hat{\mathcal{F}}_{j-1/2}^{\bmu,{\bf w}, \btheta}\right),
\end{equation}
where $\hat{\mathcal{F}}_{j+1/2}^{\bmu,{\bf w},\btheta}$ is defined in \eqref{eq: neural-entropy stable flux}.

\begin{remark}\label{remark:diffusionD}
The classical Rusanov-flux uses a scalar multiple of the identity matrix for diffusion, based on the largest eigenvalue magnitude of the Jacobian (see \eqref{eq: entropy D mat}). While simple, this choice can lead to overly diffusive behavior. A more refined alternative is the Roe-type flux, which uses the full spectral decomposition of the Jacobian of the flux function to construct the diffusion matrix:
\begin{equation} \label{eq: D_roe} 
D_{j+1/2} = \tilde{R}_{j+1/2} |\Lambda_{j+1/2}| \tilde{R}^\top_{j+1/2}, 
\end{equation} 
where $|\Lambda_{j+1/2}| = \text{diag}(|\lambda_{j+1/2}^1|, \ldots, |\lambda_{j+1/2}^p|)$ contains the absolute eigenvalues of the Jacobian. However, extending this decomposition to a learnable setting poses significant challenges, particularly in computing $\tilde{R}_{j+1/2}$ in a stable and efficient manner. We leave this direction for future investigation.
\end{remark}

\begin{remark}
Compared to the ESCFN in \cite{liu2024entropy}, which relies on a predefined classical numerical scheme, such as KT scheme (\cite{kurganov2000new}), the proposed NESCFN in this work eliminates the need for such hand-crafted baselines. The new framework is completely data-driven: all necessary information is inferred directly from the training data, including the underlying dynamics and the associated entropy-stable numerical scheme. Despite the absence of an explicit classical solver, the model retains critical structural properties through careful design of the neural architecture. 
\end{remark}

\subsection{Time integration} \label{sec: time_integrator} 

A suitable time integration method is needed to evolve the semidiscrete problem \eqref{eq: scheme_learned_law} for the learned conservation law \eqref{eq: learned_conservation_law}. Due to its reduced memory footprint, which is particularly important given the limitations of GPU memory, in our numerical experiments  we use the second-order total variation diminishing Runge--Kutta (TVDRK2) scheme \cite{Shu98}. TVDRK2 is applicable to general time-dependent ODEs as
\[\frac{d{\bf z}(t)}{dt}= \mathcal{G}\big({\bf z}(t)\big),\]
where $\mathcal{G}$ is a known operator on ${\bf z}$. Letting ${\bf z}_l$ denote the numerical approximation of ${\bf z}(t)$ at the time $t_l$,  \Cref{alg: TVDRK2} describes how TVDRK2 advances the solution from $t_{l-1}$ to $t_l = t_{l-1}+\Delta t$.


\begin{algorithm}[H]
\caption{TVDRK2 time integration method for a single time step starting at time level $t_{l-1}$}\label{alg: TVDRK2}
\begin{algorithmic}
\State INPUT: ${\bf z}_{l-1}$, $\mathcal{G}({\bf z}_{l-1})$ and $\Delta t$.
\State OUTPUT: The solution ${\bf z}_{l}$ at time level $t_l$.
\State $ {\bf z}^{(1)}={\bf z}_{l-1}+ \Delta t \, \mathcal{G}({\bf z}_{l-1})$,
\medskip
\State $ {\bf z}_{l}=\frac{1}{2}{\bf z}_{l-1} + \frac{1}{2}{\bf z}^{(1)}$.
\end{algorithmic}
\end{algorithm}

\section{Training procedure design}\label{sec: training design}
In this section, we detail the data generation process and the training protocol used in our numerical experiments presented in Section \ref{sec: numerical}. 

\subsection{Problem setup}\label{sec: data generation} Our goal is to predict the dynamics of a hyperbolic system of conservation laws by using discrete space-time solution trajectories as training data. One of the central challenges is that the flux functions governing the system's evolution are unknown. We aim to recover these dynamics from the data and use the learned model to accurately predict the evolution of conserved variables beyond the training time horizon. To clarify our presentation, we illustrate our training setup using the 1D scalar conservation law \eqref{eq:conservation_law_1d}. However, our methodology is more general, as will be demonstrated by the numerical experiments in Section~\ref{sec: numerical} for both 1D systems and the 2D Burgers equation.

\subsubsection*{Data generation and assumptions} We assume access to training data in the form of solution trajectories generated over a finite time interval. Each trajectory originates from a perturbed initial state and evolves according to the governing conservation law \eqref{eq:conservation_law_1d}. These trajectories are discretized in both time and space. In practice, such data may arise from experimental measurements or sensor networks. In this study, however, we simulate observations by numerically solving the true PDE using various perturbed initial conditions. Importantly, and in contrast to some prior work (e.g., \cite{chenPDE}), we do {\em not} rely on oracle access that might provide a richer solution space but would not be available as observations, nor do we carefully curate initial conditions to guarantee smoothness of the solution. Instead, we adopt a more realistic and general data acquisition framework, in line with the assumptions in \cite{chencfn,liu2024entropy}.

\subsubsection*{Temporal and spatial discretization} We discretize the temporal domain using a fixed time step $\Delta t$ and define by $L$ the total number of simulation steps. The full time span for data collection is then given by
\begin{equation}
    \label{eq:totalperiod}
    \mathcal{D}_{\text{train}}= [0,L\Delta t],
\end{equation} 
over which we generate trajectories from $N_{\text{traj}}$ different initial conditions. For each trajectory indexed by $k \in \{1,\cdots, N_{\text{traj}}\}$, we extract a training subinterval
\begin{equation}
    \label{eq:trainingperiod}
    \mathcal{D}^{(k)}_{\text{train}}= [t_0^{(k)},t^{(k)}_{L_{\text{train}}}].
\end{equation} 
where the starting time $t_0^{(k)}$ is sampled from the interval $[0, (L-L_{\text{train}})\Delta t]$, and the terminal time is given by $t^{(k)}_{L_{\text{train}}} = t_0^{(k)} + L_{\text{train}}\Delta t$. The observed data along each trajectory thus consists of discrete space-time samples of the conserved states
\begin{equation}
\label{eq: trajectory data}
\bu(t_l^{(k)}) \in {\mathbb R}^{n_{\text{train}}\times p}, \ l=0, \cdots, L_{\text{train}}, \ k=1, \cdots, N_{\text{traj}},
\end{equation}
where $n_{\text{train}}$ denotes the number of spatial grid points used during training, and $p$ is the number of the conserved states.


\subsection{Loss function} \label{sec:recurrentloss} 
We are now ready to present the loss function used to impose data fidelity and entropy structure on the learned conservation laws. For ease of presentation we annotate the procedure that combines \eqref{eq: scheme_learned_law} with the TVDRK2 time integrator in \cref{alg: TVDRK2} by $\mathcal N$ to represent the neural net operator for the solution update. That is, given the current net prediction from $t_{l-1}^{(k)}$, $\hat{\bu}(t_l^{(k)}) \in \mathbb{R}^{n_{\text{train}} \times p}$, the next step net prediction has the form
\begin{equation}
    \label{eq:Ndefine}
    \hat\bu(t_{l+1}^{(k)}) = \mathcal{N}\big(\hat\bu(t_{l}^{(k)})\big) \quad \text{ with } \quad \hat\bu(t_0^{(k)}) = \bu(t_0^{(k)}).
\end{equation}
 Following \cite{chencfn}, we define the \emph{recurrent loss function} as
\begin{equation}
\label{eq:recurrentloss}
\mathcal{L}\left( \btheta, \bmu, {\bf w}; \hat{\bf u} \right) = \sum_{k=1}^{N_{\text{traj}}} \sum_{l=0}^{L_{\text{train}}} \left\|\hat\bu( t_{l}^{(k)} ;\btheta, \bmu, {\bf w}) - \bu( t_{l}^{(k)} )\right\|_2^2, 
\end{equation}
 where 
 \[\hat\bu( t_{l}^{(k)};\btheta, \bmu, {\bf w} ) = \underbrace{\mathcal{N}\circ\ldots \circ \mathcal{N}}_{l \text{ times }}\big( \bu( t_{0}^{(k)}) \big).\]

A two-stage learning strategy is adopted to ensure that the learned entropy function satisfies the properties stated in Theorem~\ref{thm: existence entropy}. In the first stage we define the loss function as
\begin{equation}
\begin{aligned}
&\mathfrak{L}_1(\btheta,\bmu,{\bf w}; \hat{\bf u}) = \frac{1}{\sum_{k=1}^{N_{\text{traj}}} \sum_{l=0}^{L_{\text{train}}} \|\hat\bu( t_{l}^{(k)} )\|_2^2} \Biggl[\mathcal{L}\left( \btheta, \bmu, {\bf w};\hat{\bf u} \right)  \\
& +  \lambda_1 \sum_j\|\eta_{\btheta}''(\hat\bu_{j}) - \mathbb{I}_p\|_F^2  + \lambda_2 \sum_j\frac{\left\|\eta_{\btheta}''(\hat\bu_{j})(\mathcal{F}^{\bmu})'(\hat\bu_{j}) - [(\mathcal{F}^{\bmu})']^{\top}[(\hat\bu_{j})\eta_{\theta}'']^{\top}(\hat\bu_{j})\right\|_F^2}{{4\|\eta''_{\btheta}(\hat\bu_{j})(\mathcal{F}^{\bmu})'(\hat\bu_{j})\|_F^2}}\Biggr] ,
\end{aligned}
\label{eq: stage 1 loss}
\end{equation}
where $\mathcal{L}(\btheta, \bmu,\bf w;\hat{\bf u})$ defined in \eqref{eq:recurrentloss} enforces data fidelity. The penalty terms on the second line ensure that learned entropy function \(\eta_\btheta\) is nontrivial at the cell center $x_j$. For notational simplicity we omit the explicit dependence of $\hat {\bf u}_j$ on $t_l^{(k)}$.

The second training stage further ensures that the symmetry condition in Theorem~\ref{thm: existence entropy} is satisfied, along with the hyperbolicity of the learned conservation law. This is accomplished by fine-tuning {\em only} the final layer of the convex neural network $\eta_\btheta$ using the loss function given by
\begin{equation}
\mathfrak{L}_{2}({\btheta}_{\mathrm{L}}; \hat{\bf u}) = 
\lambda_1 \sum_j\|\eta_{\btheta}''(\hat\bu_{j}) - \mathbb{I}_p\|_F^2 
+ \sum_j\frac{\left\|\eta_{\btheta}''(\hat\bu_{j})(\mathcal{F}^{\bmu})'(\hat\bu_{j}) - [(\mathcal{F}^{\bmu})']^{\top}[(\hat\bu_{j})\eta_{\btheta}'']^{\top}(\hat\bu_{j})\right\|_F^2}{{4\|\eta''_{\btheta}(\hat\bu_{j})(\mathcal{F}^{\bmu})'(\hat\bu_{j})\|_F^2}}.
\label{eq: stage 2 loss}
\end{equation}
Here $\btheta_\mathrm{L}$ denotes the final layer trainable parameters for $\eta_\btheta$. The full two-stage training procedure for a single epoch with batch size $N_b$ is summarized in \cref{algorithm: two stage}.

\begin{algorithm}[h!]
\caption{Training Loop}
\label{algorithm: two stage}
\begin{algorithmic}[1]
\State 
\textbf{Input}: parameters \(\{\btheta, \bmu, {\bf w}\}\)
\State 
\textbf{Output}: updated parameters \(\{\btheta, \bmu, {\bf w}\}\)
\Procedure{Train in one epoch}{}
    \For{\(j = 1, \ldots, N_b\)}
    \State \(\{\btheta, \bmu, {\bf w}\}\) $\gets$ \(\{\btheta,\bmu,{\bf w}\} -\tau_1\nabla_{\btheta,\bmu, {\bf w}}\mathfrak{L}_1(\btheta, \bmu, {\bf w};\hat{\bf u})\)
        \For{$i = 1$ to stage$2_{-}$steps}
            \State  $\tilde{\bf u}\gets \hat {\bf u} + \boldsymbol{\delta}$ where $ \boldsymbol{\delta} \in \mathfrak{N}(0, 0.3\overline{|\hat\bu|})$\footnotemark 
            \State  $\btheta_{\mathrm{L}}\gets \btheta_{\mathrm{L}} - \tau_2 \nabla_{\btheta_{\mathrm{L}}} \mathfrak{L}_2\left( \btheta_{\mathrm{L}};\tilde{\bf u}\right)$ 
        \EndFor
    \EndFor
\EndProcedure
\end{algorithmic}
\end{algorithm}
\footnotetext{$\mathfrak{N}(\cdot,\cdot)$ denotes standard normal distribution, and $\overline{|\hat\bu|}$ represents the mean absolute value of $\hat{\bf u}$.}

\subsection{Conservation and entropy metric}\label{sec:conservationnmetric}
Motivated by the evaluation strategy in \cite{liu2024entropy}, we assess the performance of the proposed NESCFN from two complementary perspectives. First, we evaluate the model's prediction accuracy by comparing its forecasts with reference trajectories beyond the training horizon. Second, we verify the preservation of essential structural properties intrinsic to hyperbolic conservation laws. Specifically, we examine whether the learned model preserves the conservation of physical quantities and maintains entropy stability throughout the simulation.

To quantify conservation fidelity, we define the discrete conservation error for the $i$-th conserved variable at time $t_l$ as
\begin{equation}
\label{eq:conservemetric}
\mathcal{C}_i\big(\hat\bu(t_l)\big) := \Big| \sum_{j=1}^{n-1} \left( \hat{u}_j^i(t_l) - \hat{u}_j^i(t_0) \right)\Delta x - \sum_{s=1}^{l}\left( F_a^{i,s-1} -F_b^{i,s-1} \right)\Delta t\Big|,
\end{equation}
where $[\hat u^i_1(t_l), \ldots, \hat u^i_{n-1}(t_l)]$ denotes the predicted $i$-th conserved state at time \(t_l\), and \(F_a^{i,s-1}\) and \(F_b^{i,s-1}\) are the calculated fluxes at the domain boundaries $x = a$ and $x = b$, respectively, defined as
\begin{equation}
\label{eq:fluxterm}
        F_a^{i,s-1} = \frac{1}{\Delta t}\int_{t_{s-1}}^{t_{s}} \mathcal{F}^{\bmu,i}\big(\hat{\bf u}(a,t)\big)dt,\quad
        F_b^{i,s-1} = \frac{1}{\Delta t}\int_{t_{s-1}}^{t_{s}} \mathcal{F}^{\bmu,i}\big(\hat{\bf u}(b,t)\big)dt,
\end{equation}
where $\mathcal{F}^{\bmu,i}$ denotes the $i$-th component of the learned flux defined in \eqref{eq: learned_conservation_law}.

To assess entropy stability, we define the discrete entropy remainder at time $t_l$ as\vspace{-0.4cm}
\begin{equation}
\mathcal{J}(\hat\bu (t_l)) :=  \sum_{j=1}^{n-1}\Delta x \Big[\Big( {\eta}_{\btheta}\big(\hat{\bf u}_j(t_l)\big) - {\eta}_{\btheta}\big(\hat{\bf u}_j(t_0)\big) \Big)
- \sum_{s=1}^{l} \left( [\eta_{\btheta}^{\prime} \big(\hat{\bf u}_{j}(t_s)\big)]^\top\mathcal{F}^{\bmu}\big( \hat{\bf u}_{j}(t_s)\big) \right) \Delta t \Big] , 
    \label{eq: discrete entropy}
\end{equation}
where the second term approximates the entropy flux. We say that the network defines a {\em neural-entropy-stable} operator if $\mathcal{J}(\hat\bu (t_l)) \leq 0$. 

\subsection{Network architectures and training configuration}
\label{sec:networkdetails}

We now describe the neural network architectures used in our framework, along with the associated training parameters. As discussed in Section~\ref{sec:learning entropy stable}, the learned entropy-stable numerical flux~\eqref{eq: neural-entropy stable flux} involves three networks: $\mathcal{F}^{\bmu}$, $\eta_{\btheta}$, and $\rho_{\bf w}$. The networks $\mathcal{F}^{\bmu}$ and $\rho_{\bf w}$ are implemented as standard fully connected neural networks (FCNNs), while $\eta_{\btheta}$ is realized as an input convex neural network, a structured subclass of FCNNs designed to enforce input convexity. Both architectures are breifly reviewed below.

\subsubsection{Fully connected neural networks} A fully connected neural network (FCNN) is a feedforward architecture composed of multiple layers, where each layer applies an affine transformation followed by a nonlinear activation. Let ${\bf z}^{(0)} = \mathbf{x} \in \mathbb{R}^{d_0}$ denote the input. An $\mathrm{L}$-layer FCNN computes the output $\mathbf{y} \in \mathbb{R}^{d_\mathrm{L}}$ through the recursive relation
\begin{equation}
\begin{aligned}
\mathbf{z}^{(l)} &= \sigma\Big(\mathbf{W}^{(l)} \mathbf{z}^{(l-1)} + \mathbf{b}^{(l)}\Big), \quad l = 1,\cdots, \mathrm{L}-1,\\
{\bf y} & = \mathbf{W}^{(\mathrm{L})} \mathbf{z}^{(\mathrm{L}-1)} + {\bf b}^{(\mathrm{L})},
\end{aligned}\label{eq:activation}
\end{equation}
where $\mathbf{z}^{(0)} = \mathbf{x}$ is the input, $\mathbf{W}^{(l)} \in \mathbb{R}^{d_l \times d_{l-1}}$ and $\mathbf{b}^{(l)} \in \mathbb{R}^{d_l}$ are the weight matrices and bias vectors at layer $l$, and $\sigma(\cdot)$ is an elementwise nonlinear activation function. FCNNs are general-purpose function approximators and are used in our framework for the learned flux function in \eqref{eq: Param entropy conservative flux} and the maximum wave speed in \eqref{eq: learned_conservation_law}.


\subsubsection{Input convex neural networks}\label{sec:input_cnn}
Input convex neural networks (ICNNs) are in contrast a subclass of FCNNs designed to represent functions that are convex with respect to their input. Introduced in~\cite{amos2017input}, ICNNs impose two key structural constraints:
\begin{itemize}
\item[(1).] The activation functions must be convex and non-decreasing;
\item[(2).] The weight matrices associated with hidden-layer recursion must be element-wise non-negative to preserve convexity.
\end{itemize}
Formally, for an ICNN with input ${\bf z}^{(0)} = \mathbf{x} \in \mathbb{R}^{d_0}$ and layer outputs $\mathbf{z}^{(l)} \in\mathbb{R}^{d_l}$, the forward pass is given by
\begin{equation} 
\begin{aligned}
{\bf z}^{(l)} &= \sigma\Big({\bf W}^{z,(l)} {\bf z}^{(l-1)} + {\bf W}^{x,(l)} {\bf x} + {\bf b}^{(l)}\Big), \quad l = 1,\cdots, \mathrm{L}-1,\\
{\bf y} & = \mathbf{W}^{(\mathrm{L})} \mathbf{z}^{(\mathrm{L}-1)} + {\bf b}^{(\mathrm{L})},
\end{aligned}\label{eq: ICNN structure}
\end{equation}
where ${\bf y} \in \mathbb{R}^{d_\mathrm{L}}$ is the output, $\mathbf{W}^{z,{(l)}} \geq 0$ element-wise and $\sigma(\cdot)$ convex and non-decreasing. To maintain unconstrained optimization, the non-negativity of $\mathbf{W}^{z,(l)}$ is often enforced via an exponential reparameterization with $\mathbf{W}^{z,(l)} = e^{\tilde{\mathbf{W}}^{z,(l)}}$. ICNNs are particularly valuable when convexity is a modeling requirement, as in our case, where the entropy function $\eta_{\btheta}(\mathbf{u})$ must be convex to guarantee the hyperbolicity of the learned system.

 \subsubsection{Training configuration} We now describe the training configuration used in this work. As discussed in the preceding sections, the neural operator responsible for solution updates, denoted by $\mathcal{N}$ in~\eqref{eq:Ndefine}, is implemented using fully connected neural networks ($\mathcal{F}^\bmu, \rho_{\bf w}, \eta_{\btheta}$). Each network has both input and output dimensions equal to the number of state variables $\bf u$, i.e., $d_0 = d_{\mathrm{L}} = p$ in \eqref{eq:activation}--\eqref{eq: ICNN structure}. The architectural details of the individual networks are summarized as follows:
 \begin{itemize}
 \item[(a).] The FCNN $\mathcal{F}^\bmu$ consists of three hidden layers, each with $64$ neurons. This architecture balances expressiveness with memory constraints. We employ the Sigmoid Linear Unit (SiLU) activation function, defined as $\text{SiLU}(x) = \text{Sigmoid}(x)\cdot\text{ReLU}(x)$ in the hidden layers. The choice of SiLU is motivated by its differentiability, which is essential for computing the Jacobian of $\mathcal{F}^\bmu$ in \eqref{eq:activation};
 \item[(b).] The FCNN $\rho_{\bf w}$ uses a simpler architecture with two hidden layers of $64$ neurons each. Since higher-order differentiability is not required in this case, we adopt the ReLU activation function $\text{ReLU}(x) = \max\{0, x\}$ for its simplicity and computational efficiency;
\item[(c).] The ICNN $\eta_\btheta$ is modeled using a one-layer ICNN with $64$ hidden neurons. To ensure convexity and smoothness, we apply the Softplus activation function, defined as $\text{Softplus}(x) = \log(1+e^x)$, which is convex, non-decreasing and continuously differentiable.
 \end{itemize}

The Adam optimizer (\cite{Adam}) is employed during training to update the neural networks' weights and biases. To maintain consistency across experiments, we use fixed learning rates $\tau_1$ and $\tau_2$ for the two training stages in \Cref{algorithm: two stage}, without employing any learning rate scheduler. \Cref{tab:training parameter} provides the full set of training hyperparameters  for each experiment conducted in \Cref{sec: numerical}: the loss weights $\lambda_1$ and $\lambda_2$ in \cref{eq: stage 1 loss} and \cref{eq: stage 2 loss}, the number of training epochs $N_{\text{Epoch}}$, the batch size $N_b$, the learning rates $\tau_1$ and $\tau_2$ for the two training stages in \Cref{algorithm: two stage}, and the number of fine-tuning steps used in stage 2. 
\begin{table}[h!]
\centering
\begin{tabular}{|l|c|c|c|c|c|c|c|}
\hline
 & $\lambda_1$ & $\lambda_2$ & \small{$N_{\text{Epoch}}$} & \small{$N_b$} & $\tau_1$ & \small{stage2\_steps} & $\tau_2$ \\
\hline
\small{1D Burgers' eq}  & $10^{-3}$  & $10^{-2}$ &  $50$ & $5$  & $1\times 10^{-3}$  & $0$  & --  \\
\hline
\small{1D shallow water eq}  & $10^{-3}$  & $10^{-2}$  & $100$  &  $5$ & $2\times 10^{-3}$  & $0$  &  -- \\
\hline
\small{1D Euler's eq}  &  $10^{-5}$ & $10^{-3}$  & $500$  &  $10$ & $2\times 10^{-3}$ & $5$   & $10^{-3}$  \\
\hline
\small{2D Burgers' eq}   & $10^{-5}$  & $10^{-2}$  & $50$  & $1$  & $1\times 10^{-3}$  & $0$  &  -- \\
\hline
\end{tabular}
\caption{Training parameters for each experiment conducted in \Cref{sec: numerical}.}
\label{tab:training parameter}
\end{table}

The choice of hyperparameters in Table~\ref{tab:training parameter} reflects a balance between the complexity of the problem, training efficiency, and available GPU memory. While the number of training trajectories, \(N_{\text{traj}}\), varies across problems, we fix the number of validation cases to $40$ in all experiments, which is used to determine whether the model parameters should be saved at each epoch. Additionally, the constants \(C_1 = 10\) and \(C_2 = 0.3\) in \eqref{eq: stable cho decomp} are fixed for all cases to promote rapid decay. All implementations are based on JAX (\cite{jax}), a high-performance numerical computing library that supports automatic differentiation and GPU/TPU acceleration in Python. Finally, to ensure robustness, we note that none of these parameters were fine-tuned.\footnote{The complete code  is  available upon request for reproducibility purposes.}

\section{Numerical experiments}
\label{sec: numerical}

We now provide a series of numerical experiments to demonstrate that the NESCFN offers a flexible procedure to learn the \textit{entropy} and the corresponding \textit{entropy-stable} scheme. In particular, we show that the NESCFN can predict not only the long-term behavior of the dynamics, but that it also preserves the entropy inequality for all noise levels in our experimental training data. We show this for both in-distribution initial conditions and out-distribution initial conditions, highlighting the robustness and generalizability of the proposed learning framework. Specifically, the out-distribution initial conditions considered here include two scenarios. The first involves initial conditions that share a similar structural pattern with the training data, but whose coefficients are sampled from a distribution different from that used during training (see e.g., \eqref{eq: Burgers' initial condition} and \eqref{eq:entropy_test_1dburger}). The second involves initial conditions with qualitatively different patterns from those in the training set -- for example, testing on Gaussian profiles when the training data were generated using sinusoidal waves (see  e.g., \eqref{eq: Burgers' initial condition 2d} and \eqref{eq:entropy_test_2dburger2}). Both types of out-distribution data are included in our evaluation of neural entropy stability.



\subsection{Prototype conservation laws}
\label{sec: pdeexamples}
The classical conservation laws described below  serve as prototype test cases to evaluate the effectiveness and robustness of the proposed NESCFN.


\subsubsection{1D Burgers' equation}
\label{sub:1DBurgers}
We first consider the scalar Burgers' equation
\begin{equation}
    \dfrac{\partial u}{\partial t} + \dfrac{\partial}{\partial x} \Big( \dfrac{u^2}{2} \Big) = 0, \quad x \in \left[ 0, 2\pi \right], \quad t > 0,
    \label{eq: Burgers' equation}
\end{equation}
with periodic boundary conditions $u(0,t)=u(2 \pi, t)$, and initial condition
\begin{equation}
    u\left( x, 0 \right) = \alpha \sin\left( x \right) + \beta, \quad \alpha,\beta \in \mathbb{R}.
    \label{eq: Burgers' initial condition}
\end{equation}
where the parameters in \eqref{eq: Burgers' initial condition} are sampled uniformly as $\alpha \sim \mathcal{U}\left[ .75, 1.25 \right]$ and $\beta \sim \mathcal{U}\left[ -.25, .25 \right]$ to generate training data used in our experiments.

As discussed in \Cref{sec: data generation}, each of the $N_{\text{traj}}$ training trajectories is generated using the PyClaw package (\cite{clawpack, pyclaw}), with a fixed time step \(\Delta t = .005\) and spatial grid size $n_{\text{train}} = 512$ . We set $N_{\text{traj}}=200$, and define the total training time period (see \cref{eq:totalperiod}) as $L = 20$. We do not subdivide the training interval into smaller temporal segments in this example; that is, we take $L_{\text{train}} = L =20$ in \cref{eq:trainingperiod}. Observe that under these settings only \emph{smooth} solution profiles prior to shock formation are used for training.

For testing, we fix $\alpha = 1.05609$ and $\beta = 0.1997$ in \cref{eq: Burgers' initial condition}. The corresponding reference solution is computed using PyClaw with the same temporal and spatial resolution as used in the training data generation, namely, $(\Delta t, n_{\text{test}}) = (.005, 512)$, and integrated up to time $T = 3$.

To assess learned entropy stability, we 
consider a family of initial conditions
\begin{equation}\label{eq:entropy_test_1dburger}
    u\left( x,0 \right) = (0.5 + 0.01\kappa)\sin \bigg[ \Big(1 + \Big\lfloor \frac{\kappa}{20} \Big\rfloor\Big) x + 0.01\kappa \bigg],
\end{equation}
where \(\kappa = 0,1,\ldots,N_{\text{ent}}\), and \(N_{\text{ent}} = 100\) denotes the total number of trajectories used for entropy stability evaluation. This set of initial conditions includes so called out-distribution examples relative to the training data and reflects practical scenarios where data may originate from heterogeneous measurement sources with conflicting or non-overlapping characteristics. Evaluating model performance under such conditions provides valuable insight into its generalizability and robustness.

\subsubsection{Shallow water equation} \label{sub:shallowwater} We next consider the 1D shallow water system defined over the spatial domain $x\in(-5,5)$ given by 
\begin{equation}
    \begin{aligned}
        \dfrac{\partial h}{\partial t} + \dfrac{\partial}{\partial x} \left( hu \right) &= 0, \\
        \dfrac{\partial}{\partial t} \left( hu \right) + \dfrac{\partial}{\partial x} \Big( hu^2 + \dfrac{1}{2}gh^2 \Big) &= 0,
    \end{aligned}
    \label{eq: shallow water equation}
\end{equation}
subject to Dirichlet boundary conditions. The initial conditions follow
\begin{equation}
    h\left( x, 0 \right) = \left\{
    \begin{aligned}
     &h_l + \omega_{l}, \quad x < x_0 + \omega_{x}   \\
     &h_r + \omega_{r}, \quad x \geq x_0 + \omega_{x} \\
     \end{aligned}
    \right.,\quad 
    u\left( x, 0 \right) = \left\{
    \begin{aligned}
     &u_l+\omega_{ul}, \quad x < x_0+\omega_{x}   \\
     &u_r+\omega_{ur}, \quad x \geq x_0+\omega_{x} \\
     \end{aligned}
    \right..\label{eq: shallow water condition}
\end{equation}
Here, \(h_l = 3.5, h_r=1.0, u_l=u_r=x_0=0\), and \(\omega_{l}, \omega_{r}, \omega_{ul}, \omega_{ur}, \omega_{x} \in \mathbb{R} \).

For training, the parameters in \eqref{eq: shallow water condition} are sampled uniformly as \(\omega_{l}, \omega_{r}\sim \mathcal{U}\left[ -.2, .2 \right]\) and \(\omega_{ul}, \omega_{ur}, \omega_{x} \sim \mathcal{U}\left[ -.1,.1 \right]\). All other training parameters are kept consistent with those used in Section~\ref{sub:1DBurgers} with $(\Delta t, n_{\text{train}}, N_{\text{traj}}, L, L_{\text{train}}) = (.005,512, 300, 20, 20)$. Each trajectory is computed using PyClaw with the HLLE Riemann solver.

For testing, we fix the parameters in \eqref{eq: shallow water condition} to $h_l = 3.5691196, h_r = 1.178673$, $u_l = -.064667, u_r = -.045197$, and $x_0 = .003832$ with \(\omega_{l}=\omega_{r}=\omega_{ul}=\omega_{ur}=\omega_{x}=0\). 
The corresponding reference solution is obtained using the same temporal and spatial resolution as in the training data generation and is simulated up to time $T = 1.5$. 

As before we consider the case where conflicting training data may be observed, modeled here using out-distribution initial conditions. It is important to ensure that the NESCFN is able to predict long-term behavior correctly and preserve the entropy inequality in this type of training environment. To this end, we evaluate the model using the following family of initial conditions:
\begin{equation}\label{eq:entropy_test_1dshallow}
    h\left( x, 0 \right) = \begin{cases}
     6.0 - 0.01\kappa, \quad x < 0 + 0.01\kappa,   \\
     0.1 + 0.01\kappa, \quad x \geq 0 + 0.01\kappa,       
    \end{cases} 
    \quad u\left( x, 0 \right) = 0,
\end{equation}
with \(\kappa = 0,1,\ldots,N_{\text{ent}} = 100\).



\subsubsection{Euler's equation}\label{sub:Eulers} We now consider the system of 1D Euler's equations
\begin{equation}
    \begin{aligned}
        \rho_{t} + \left( \rho u \right)_{x} &= 0, \\
        \left( \rho u \right)_{t} + \Big( \rho u^2 + p \Big)_{x} &= 0, \\
        \left( E \right)_{t} + \left( u (E + p) \right)_{x} &= 0,
    \end{aligned}
    \label{eq: Euler's equation}
\end{equation}
for $x \in (-5,5)$ with Dirichlet boundary conditions. To  demonstrate its robustness across varying levels of complexity, we will evaluate our NESCFN framework on two classic initial conditions, the Sod shock tube and the more challenging Shu–Osher problem. Both are designed to include out-distribution components relative to the training data, thereby simulating realistic scenarios with conflicting sources of observational data. In this regard, the general training initial condition that we will use for in-distribution testing is defined as
\begin{equation}
\begin{aligned}
& \rho(x, 0)= \begin{cases}\rho_l, & \text { if } x \leq x_0, \\
1+\varepsilon \sin (5 x), & \text { if } x_0<x \leq x_1, \quad u(x, 0)= \begin{cases}u_l, & \text { if } x \leq x_0, \\
0, & \text { otherwise, }\end{cases} \\
1+\varepsilon \sin (5 x) e^{-\left(x-x_1\right)^4}, & \text { otherwise, }\end{cases} \\ \\
& p(x, 0)=\left\{
\begin{array}{ll}
p_l, & \text { if } x \leq x_0, \\
p_r, & \text { otherwise, }
\end{array} \quad E(x, 0)=\frac{p_0}{\gamma-1}+\frac{1}{2} \rho(x, 0) u(x, 0)^2, \right.
\end{aligned}\label{eq:euler IC}
\end{equation}
where $x_1=3.29867$ and $\gamma=1.4$. 

To generate training trajectories, we sample the parameters in \eqref{eq:euler IC} uniformly as\vspace{-0.2cm}
\[\begin{+array}{lrlr}
 \rho_l \sim \mathcal{U}[\hat{\rho}_l(1-\epsilon), &\hat{\rho}_l(1+\epsilon)], &
 \varepsilon\; \sim \mathcal{U}[\hat{\varepsilon}(1-\epsilon),  &\hat{\varepsilon}(1+\epsilon)], \\
 p_l \sim \mathcal{U}[\hat{p}_l(1-\epsilon), &\hat{p}_l(1+\epsilon)], &
 p_r \sim \mathcal{U}[\hat{p}_r(1-\epsilon),  &\hat{p}_r(1+\epsilon)], \\
 u_l \sim \mathcal{U}[\hat{u}_l(1-\epsilon), &\hat{u}_l(1+\epsilon)], & 
 x_0 \sim \mathcal{U}[\hat{x}_0(1-\epsilon), &\hat{x}_0(1+\epsilon)],
\end{+array}\vspace{-0.2cm}\]
with $\epsilon=.1$, $\hat{\rho}_l = 3.857135$, $\hat{p}=10.32333$, $\hat{u}_l=2.62936$, $\hat{\varepsilon} =.2$, $\hat{p}_r=1$, and $\hat{x}_0=-4$. Training trajectories are computed using PyClaw with HLLE Riemann solver, a fixed time step \(\Delta t = .002\), spatial resolution \(n_{\text{train}} = 512\), and a total of $N_{\text{train}} = 300$ samples. The total training time domain $\mathcal{D}_{\text{train}}$ (see \eqref{eq:totalperiod}) is fixed with \(L = 300\). To enrich the training dataset and expose the model to a broader range of local dynamics, we partition the full time domain $\mathcal{D}_{\text{train}}$ into shorter overlapping segments. Specifically, we set $L_{\text{train}} = 20$ in \eqref{eq:trainingperiod} and define the starting time of each segment as $t_0^{(k)} = l\Delta t$, where $l \in \{0,1, \cdots, L-L_{\text{train}}\}$. This approach yields $N_{\text{train}} \times (L- L_{\text{train}} + 1)$ training trajectories, significantly more than the $N_{\text{train}}$ trajectories used in the Burgers' and shallow water equation examples, where the full training window is used with $L_{\text{train}} = L$. The increased number of shorter trajectories enhances the training diversity while reducing the risk of overfitting to long-term profiles. 

For testing we set \(\rho_l = \hat{\rho}_l, \varepsilon = \hat{\varepsilon}, p_l = \hat{p}_l, p_r = \hat{p}_r, u_l = \hat{u}_l,\) and \(x_0 = \hat{x}_0\) in \eqref{eq:euler IC}, and evolve the system \eqref{eq: Euler's equation} up to \(T = 1.6\) using $\Delta t = .002$ and \(n_{\text{test}} = 512\). 

To validate our new NESCFN  approach, we test on both the Sod shock tube and Shu-Osher problems (\cite{SHU1988439}), each designed to include out-distribution initial conditions relative to the training set. 

The initial conditions for the Sod problem are given by
\begin{equation}\label{eq:entropy_test_1deuler_rod}
\begin{array}{lll}
&\!\! \rho(x, 0)\!\!= \!\!\begin{cases}
3.5, & \text{ if } x \leq x_0, \\
0.12 + 0.01\kappa, & \text{ otherwise,}
\end{cases} 
&\!\! p(x, 0)\!\!= \!\!\begin{cases}
10.0 - 0.01\kappa, & \text{ if } x \leq x_0, \\
1.0 + 0.01\kappa, & \text{ otherwise,}
\end{cases} \\ \\
&\!\! E(x, 0)\!\!=\!\!\frac{p_0}{\gamma-1}+\frac{1}{2} \rho(x, 0) u^2(x, 0), 
&\!\! u(x, 0) \!\!=\! 0,
\end{array}
\end{equation}
with $\kappa = 0,1,\cdots,N_{\text{ent}}$, and the interface location
\(
x_0 = -0.3 + 0.6 (2\kappa)/N_{\text{ent}}.
\) We emphasize that the initial condition in \eqref{eq:entropy_test_1deuler_rod} differs structurally from the training form \eqref{eq:euler IC}, introducing sharp discontinuities and spatial variations absent in the training data. This design tests the robustness and generalizability of NESCFN to out-distribution settings.

For the Shu-Osher problem we parameterize the initial conditions \cref{eq:euler IC} with
\begin{equation}\label{eq:entropy_test_1d_shuO_osher}
\begin{array}{llll}
& \rho_l =3.857135, 
& \varepsilon = .1 + 0.005\kappa, 
& p_l = 10.33333 - 0.01\kappa, \\
& p_r = 1 + 0.01\kappa, 
& u_l = 2.629 - 0.01\kappa, 
& x_0 = -0.8 + 0.01\kappa, 
\end{array}
\end{equation}
for \(\kappa = 0, 1, \ldots, N_{\text{ent}}\).  

Here, again, \(N_{\text{ent}} = 100\) denotes the total number of trajectories used for entropy stability evaluation for both Sod and Shu-Osher problems.

\subsubsection{2D Burgers' equation}\label{sub: 2D Burgers} Finally, to evaluate the proposed NESCFN in a multidimensional setting, we consider the 2D Burgers' equation
\begin{equation}
    \dfrac{\partial u}{\partial t} + \dfrac{\partial}{\partial x} \Big( \dfrac{u^2}{2} \Big) + \dfrac{\partial}{\partial y} \Big( \dfrac{u^2}{2} \Big) = 0, \quad (x, y) \in \left[ 0, 1 \right]\times\left[ 0, 1 \right], \quad t > 0, 
    \label{eq: Burgers' equation 2d}
\end{equation}
subject to periodic boundary conditions. The initial condition is given by
\begin{equation}   
u\left( x, y, 0 \right) = \alpha \sin\left( 2\pi x +x_0\right)\cos\left( 2\pi y +y_0\right) + \beta,\quad \alpha, \beta,x_0,y_0 \in \mathbb{R}.
    \label{eq: Burgers' initial condition 2d}
\end{equation}

For training the parameters in \eqref{eq: Burgers' initial condition 2d} are sampled uniformly with \(\alpha \sim {\mathcal U}[.75,1.25]\), \(\beta \sim {\mathcal U}[-.25,.25]\),  \(x_0 \sim {\mathcal U}[.5,1.5]\), and \(y_0 \sim {\mathcal U}[-.5, .5].\) Each of the $N_{\text{traj}}$ training trajectories is again computed using PyClaw with the HLLE Riemann solver, with parameters set as $(\Delta t, n_{\text{train}}, N_{\text{traj}}, L, L_{\text{train}}) = (.001, 100\times 100, 5, 20, 20)$.

For testing we fix the parameters in \eqref{eq: Burgers' initial condition 2d} to \(x_0 = 1.032833\), \(y_0 = .034137\), \(\alpha = 1.004777\), and \(\beta = .106782.\) The reference solution is computed on a spatial grid of $n_{\text{test}} = 100\times 100$ using time step $\Delta t = .001$ integrated up to time \( T = 1.6\).

 To assess entropy stability and the model's ability to generalize beyond the training distribution, we evaluate NESCFN on three families of out-distribution initial conditions parameterized by $\kappa = 0,1,\cdots,N_{\text{ent}} = 100$. These variations simulate potential mismatches or diversity in observational sources. These families can be described as
 \begin{enumerate}
     \item  modulated sine waves 
\begin{equation}\label{eq:entropy_test_2d_burger1}
u\left( x,y,0 \right)\! =\! A\sin\left( 2 \pi k_1 x +  0.1\kappa \right) \sin\big( 2 \pi k_2 y + 0.2\kappa \big),
\end{equation}
where \(k_1 = \kappa\bmod 5 \), \(k_2 = ( 2 + \kappa )\bmod 5\), and $A = 0.5 + 0.05\kappa$. 
\item  parameterized Gaussian bumps
\begin{equation}\label{eq:entropy_test_2dburger2}
u\left( x,y,0 \right) = (0.8 + 0.02\kappa) \exp\big( - \left( x- x_0  \right)^2/(2 \sigma_{x}^2) - \left( y - y_0 \right)^2/(2 \sigma_{y}^2) \big),
\end{equation}
with \(x_0 = 0.2\left( \kappa \bmod 5 - 2  \right)\), \(y_0 = 0.2\big( \lfloor \frac{\kappa}{5} \rfloor - 2\big)\), \(\sigma_{x} = 0.1 + 0.02\left( \kappa\bmod 3 \right)\), and \(\sigma_{y} = 0.1 + 0.02 \lfloor \frac{\kappa}{3} \rfloor \). 
\item  asymmetric localized waves
\begin{equation}\label{eq:entropy_test_2dburger3}
u\left( x,y,0 \right) = (0.5 + 0.05\kappa) \left( x - x_0 \right)\exp\big( - \left( x - x_0 \right)^2 - \left( y - y_0 \right)^2  \big),
\end{equation}
where \(x_0 = 0.3\left( \kappa \bmod 3 - 1  \right)\) and \(y_0 = 0.3\big( \lfloor \frac{\kappa}{3} \rfloor -1 \big).\)
\end{enumerate}

\subsection{Observations with varying noise coefficients}
\label{subsec: noisy}
We now consider a training environment corrupted by noise, which is simulated  by applying additive noise to the  training trajectories corresponding to the PDE models introduced in \Cref{sec: pdeexamples}. 

\subsubsection{1D Burgers' equation} 
\label{subsubsec: 1D Burgers Noise}
Following the setup described in \Cref{sub:1DBurgers}, we generate  perturbed training data in domain $\mathcal{D}^{(k)}_{\text{train}}$ (see \cref{eq:trainingperiod}) as 
\begin{equation}
\tilde{u}\big( x_i, t_l^{(k)} \big) = u\big( x_i, t_l^{(k)} \big) + \xi \overline{|u\left( x,t \right)|} \zeta_{i,l}, \quad k = 1,\dots,N_{\text{traj}}.
\label{eq: noise burgers}
\end{equation}
Here \(\xi \in [0,1]\) is defined as the noise coefficient controlling the intensity of the perturbation while $\zeta_{i,l} \sim \mathfrak{N}\left( 0, 1 \right)$ denotes independent standard normal random variables, with indices $i = 1,\dots, n_{\text{train}}$ and $l = 0,\dots, L_{\text{train}}$.\footnote{Recall in this example we take $L_{\text{train}} = L$ (see \cref{eq:totalperiod} for the definition of $L$).} The quantity \(\overline{|u\left( x,t \right)|}\) represents the mean absolute value of exact solution \(u\left( x,t \right)\) over the entire dataset.  

\begin{figure}[h!]
        \includegraphics[width=\textwidth]{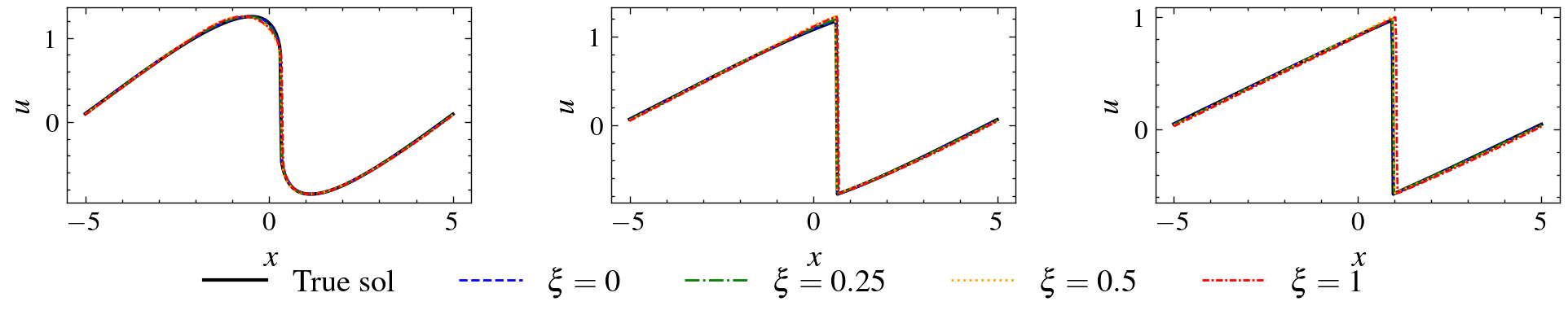}
    \caption{Comparison of the reference solution (black solid line) to 1D Burgers' equation with predictions from the NESCFN at (left) \(t = 1\) (middle) \(t = 2\) (right) \(t = 3\) for $\xi = 0,.25,.5,1$ in  \eqref{eq: noise burgers}.}
    \label{fig: 1D Burgers Predictions Noise}
\end{figure}

\Cref{fig: 1D Burgers Predictions Noise} displays a temporal sequence of the solution for the NESCFN model described in \Cref{sec:learning entropy stable} applied to the 1D Burgers' equation with noisy training data as defined in \eqref{eq: noise burgers} for $\xi = 0,.25, .5$ and $1$. Notably, the solution of the NESCFN model captures the shock formation occurring at \( t = 1.0\), {\em even though} the training data are restricted to $\mathcal{D}_{\text{train}} = [0, .1]$, a time interval prior to shock formation during which the solution remains smooth. This demonstrates that the NESCFN method does not rely on oracle knowledge of later-time dynamics and can generalize beyond the training window. Furthermore, the shock structure is preserved in long-time predictions for noise levels up to $\xi\leq0.5$, in contrast to the KT-ESCFN method proposed in \cite{liu2024entropy}, where stability is maintained even for $\xi\leq1.0$. While this comparison suggests that the NESCFN is more sensitive to structural uncertainty, we note that the entropy function is both known and fixed for the KT-ESCFN method, while here the entropy function must be learned along with the flux. In particular the large set of admissible choices for the diffusion matrix $D$ used in constructing entropy-stable fluxes introduces significant flexibility and hence additional uncertainty into the model architecture, thereby increasing the difficulty of the learning task. We defer the investigation of more sophisticated designs for $D$ to future work (see \Cref{remark:diffusionD}).

\begin{figure}[h!]
        \centering
        \includegraphics[width=0.7\textwidth]{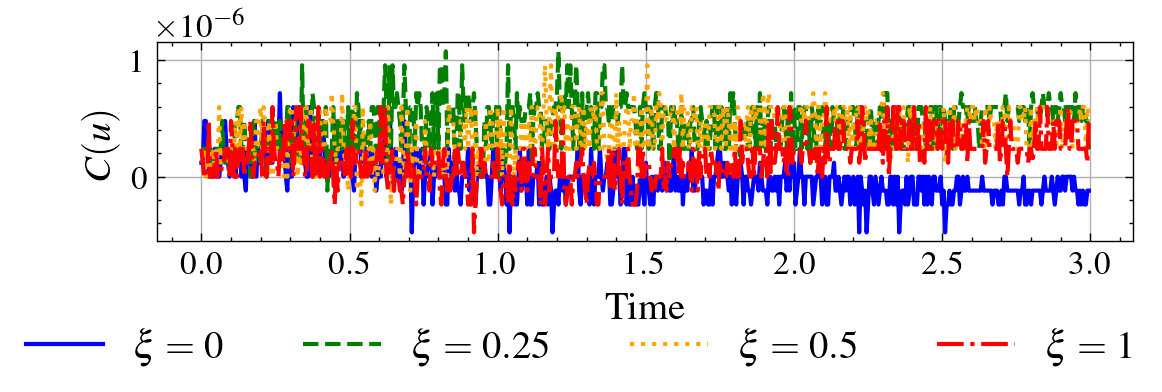}
       \caption{Discrete conserved quantity remainder $\mathcal{C}(u)$, defined in \eqref{eq:conservemetric}, for the 1D Burgers' equation over the time interval $t \in [0,3]$ with noise levels $\xi = 0, .25, .5, 1$ in  \eqref{eq: noise burgers}.} 
    \label{fig: 1D Burgers Conservations Noise}
\end{figure}

\begin{figure}[h!]
    \centering
    \begin{subfigure}{0.23\textwidth}
    \centering
        \includegraphics[width=\textwidth]{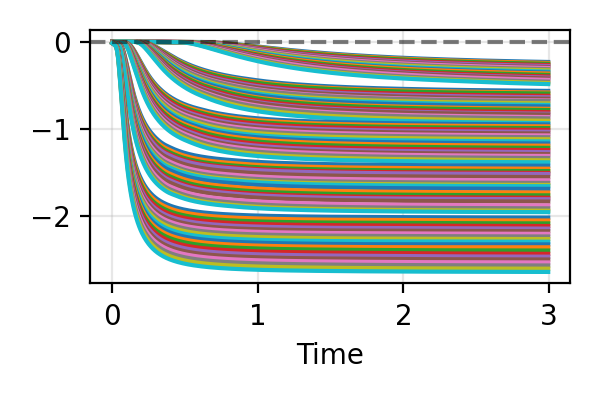}
    \end{subfigure}\hspace{0.1cm}
    \begin{subfigure}{0.23\textwidth}
    \centering
        \includegraphics[width=\textwidth]{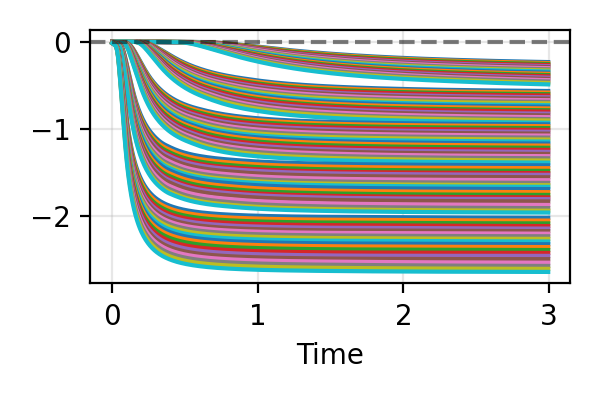}
    \end{subfigure}\hspace{0.1cm}
    \begin{subfigure}{0.23\textwidth}
    \centering
        \includegraphics[width=\textwidth]{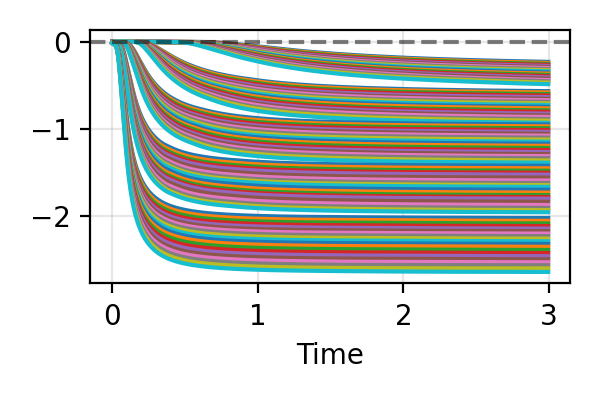}
    \end{subfigure}
    \begin{subfigure}{0.23\textwidth}
    \centering
        \includegraphics[width=\textwidth]{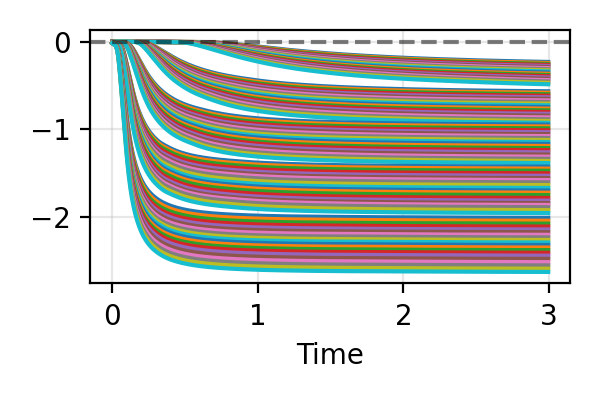}
    \end{subfigure}
    \caption{ Discrete entropy remainder $\mathcal{J}(u)$, defined in \eqref{eq: discrete entropy}, for the 1D Burgers' equation over the time interval $t \in [0,3]$, evaluated across all testing initial conditions specified in \eqref{eq:entropy_test_1dburger}. From (left) to (right) the plots correspond to training data with respective noise levels $\xi = 0, .25, .5$, and $1$ in \cref{eq: noise burgers}.} 
    \label{fig: 1D Burgers Entropy Remainder}
\end{figure}

\Cref{fig: 1D Burgers Conservations Noise,fig: 1D Burgers Entropy Remainder} further illustrate the robustness of the NESCFN method to noise.  \Cref{fig: 1D Burgers Conservations Noise} shows how the discrete conservation remainder $\mathcal{C}({u})$ as defined in \eqref{eq:conservemetric} evolves in time for different noise levels, while  \Cref{fig: 1D Burgers Entropy Remainder} demonstrates the non-positivity of the learned neural entropy $\eta_\btheta$ by displaying the corresponding evolution of the discrete entropy remainder, $\mathcal{J}({u})$ in \eqref{eq: discrete entropy}, across all testing initial conditions specified in \eqref{eq:entropy_test_1dburger}. Taken together, these results confirm that the NESCFN method provides consistent and robust long-term predictions for the 1D Burgers' equation, even in noisy training data environments.

\subsubsection{Shallow water equation} 
\label{subsubsec: Shallow Water Noise}

We now examine the impact of noise in the training data for the shallow water equations \cref{eq: shallow water equation}. As in the 1D Burgers’ case, zero-mean Gaussian noise is added to the training data within the domain $\mathcal{D}^{(k)}_{\text{train}}$, defined in \cref{eq:trainingperiod} for $k = 1,\dots,N_{\text{traj}}$. The pertubed data are given by
\begin{equation}
    \tilde{\bf a}(x_i,t_l^{(k)})
     = {\bf a}(x_i,t_l^{(k)})
     + \xi\overline{|\bm{a}|}
       \zeta_{i,l},
       \label{eq:noise_shallow}
\end{equation}
where \({\bf a} = [h,hu]^\top\) denotes the vector of physical variables, \(\zeta_{i,l} \sim \mathfrak{N}({\bf 0}, \mathbb{I}_2)\) is a 2-dimensional standard normal vector, \( i = 1,\dots, n_{\text{train}}, l = 0,\dots, L_{\text{train}}\), and \(\overline{|\bm{a}|}\) is the mean absolute value of the training data over the entire dataset.  We consider  noise intensity coefficient  \(\xi = 0, .25, .5,\) and $1$. Recall that $L_{\text{train}} = L$ in this example.

\begin{figure}[h!]
        \includegraphics[width=\textwidth]{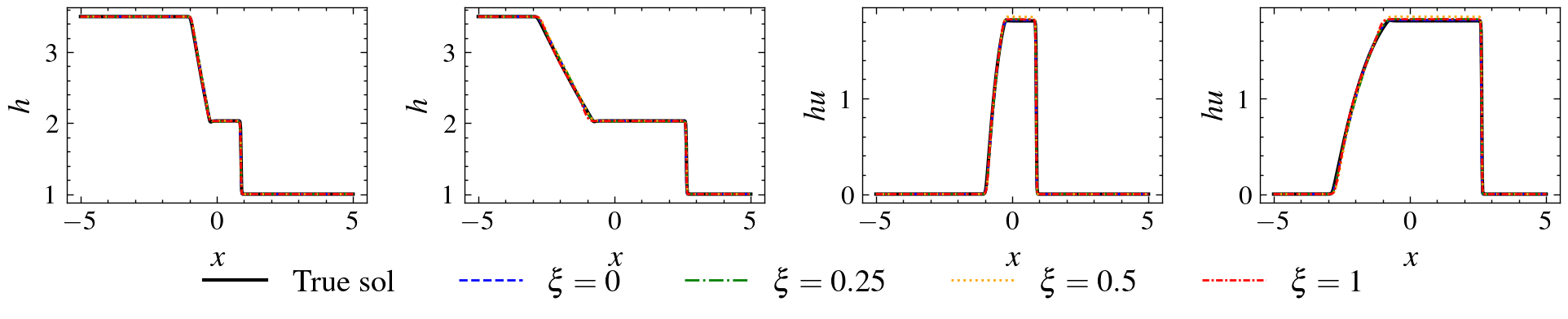}
        \caption{Comparison of the reference solution (black solid line) of height \(h\) and momentum \(hu\) for the shallow water equations with the NESCFN model predictions for  \(\xi = 0, .25, .5, 1\) in  \eqref{eq:noise_shallow}: (left) \(t = .5\) of \(h\), (middle-left) \( t= 1.5\) of \(h\), (middle-right) \(t = .5\) of \(hu\), (right) \(t = 1.5\) of \(hu\).} 
    \label{fig: Shallow Water Predictions Noise}
\end{figure}

\Cref{fig: Shallow Water Predictions Noise} displays the predicted height \(h\) and momentum \(hu\) at times \(t = .5,\) and \(t = 1.5\) for each value of $\xi$.  Remarkably, the NESCFN method accurately captures the shock structure even at the noise level $\xi = 1$. Compared to the 1D Burgers' case, the improved performance  can be attributed to the presence of discontinuities in the training data, which provides the network with direct exposure to shock-like features during training.

\begin{figure}[h!]
        \centering
        \includegraphics[width=0.6\textwidth]{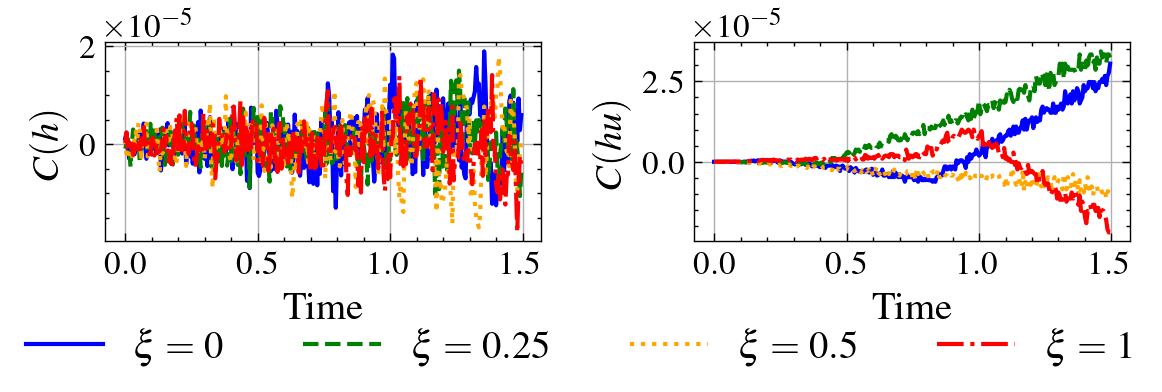}
    \caption{Discrete conserved quantity remainder  \eqref{eq:conservemetric} of height $\mathcal{C}(h)$ (left) and the momentum \(\mathcal{C}(hu)\) (right) for the shallow water equations with noise levels $\xi = 0,.25,.5, 1$ in  \eqref{eq:noise_shallow}.} 
    \label{fig: Shallow Water Conservations Noise}
\end{figure}

\begin{figure}[h!]
    \centering
    \begin{subfigure}{0.23\textwidth}
    \centering
        \includegraphics[width=\textwidth]{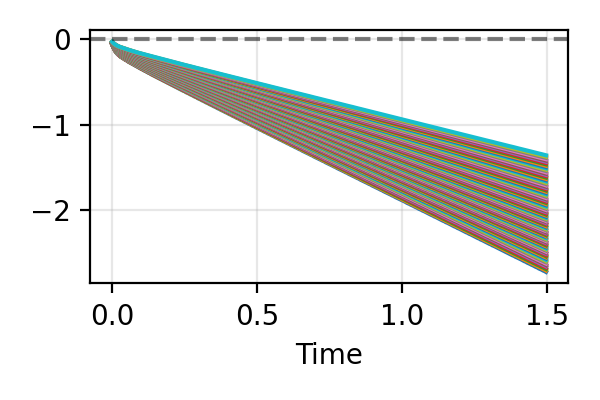}
    \end{subfigure}\hspace{0.1cm}
    \begin{subfigure}{0.23\textwidth}
    \centering
        \includegraphics[width=\textwidth]{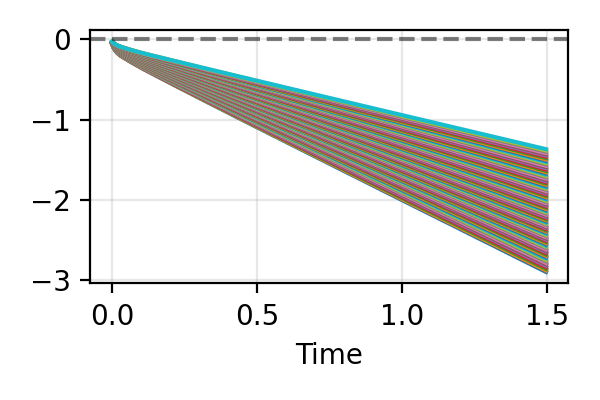}
    \end{subfigure}\hspace{0.1cm}
    \begin{subfigure}{0.23\textwidth}
    \centering
        \includegraphics[width=\textwidth]{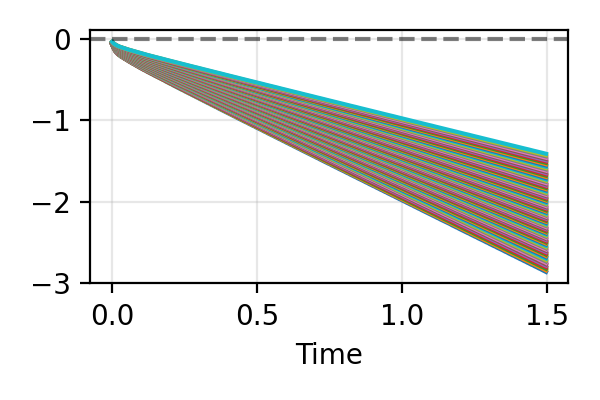}
    \end{subfigure}
    \begin{subfigure}{0.23\textwidth}
    \centering
        \includegraphics[width=\textwidth]{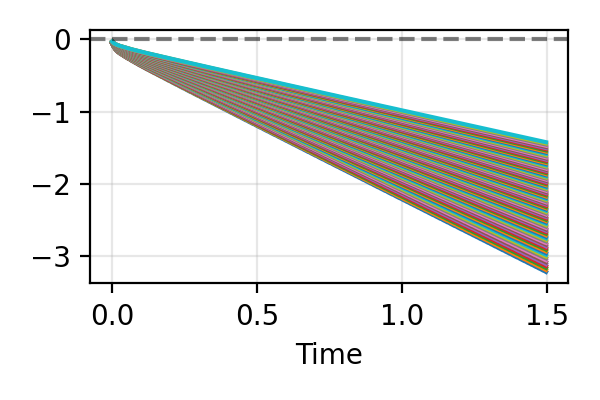}
    \end{subfigure}
    \caption{ Discrete entropy remainder \(\mathcal{J}([h,hu]^\top)\), defined in \eqref{eq: discrete entropy}, for the shallow water equations, evaluated across all testing initial conditions specified in \eqref{eq:entropy_test_1dshallow}. From (left) to (right), the plots correspond to training data with noise levels $\xi = 0, .25, .5$, and $1$, respectively.}
    \label{fig: Shallow Water Entropy Remainder}
\end{figure}

\Cref{fig: Shallow Water Conservations Noise} presents the evolution of the discrete conserved quantity remainder $\mathcal{C}(h)$ and \(\mathcal{C}(hu)\), as defined in \cref{eq:conservemetric}, while \Cref{fig: Shallow Water Entropy Remainder} shows the corresponding discrete entropy remainder $\mathcal{J}([h,hu]^\top)$ from \cref{eq: discrete entropy} across all testing initial conditions given in \eqref{eq:entropy_test_1dshallow}. The results indicate that conservation is maintained up $\mathcal{O}(10^{-5})$, and the entropy remainder remains non-positive, confirming the entropy stability of the proposed NESCFN method.

\subsubsection{Euler's equation} 
\label{subsubsec: Euler Noise}
We now investigate the impact of noise in the training data for the Euler equations \cref{eq: Euler's equation}. As before, zero-mean Gaussian noise is added to the training data within the domain \(\mathcal{D}_{\text{train}}^{(k)}\) for \(k = 1,\ldots, N_{\text{traj}}\), resulting in
\begin{equation}
    \tilde{\bf a}(x_i,t_l^{(k)})
     = {\bf a}(x_i,t_l^{(k)})
     + \xi\overline{|\bm{a}|}
       \zeta_{i,l},
       \label{eq:noise_euler}
\end{equation}
where \(\bm{a} = \left[ \rho, \rho u, E \right]^{T}\) represents the state vector, \(\zeta_{i,l} \sim \mathfrak{N}({\bf 0}, \mathbb{I}_3)\) is a standard 3-dimensional Gaussian vector, \(i = 1, \ldots, n_{\text{train}}\),  \(l = 0, \ldots,L_{\text{train}}\), with \(L_{\text{train}} < L\) (see Section \ref{sub:Eulers}). The term \(\overline{|\bm{a}|}\) denotes the mean absolute value of \({\bm a}\) over the entire dataset. We consider  noise intensity coefficient  \(\xi = 0, .25, .5,\) and $1$. 
\begin{figure}[h!]
    \includegraphics[width=\textwidth]{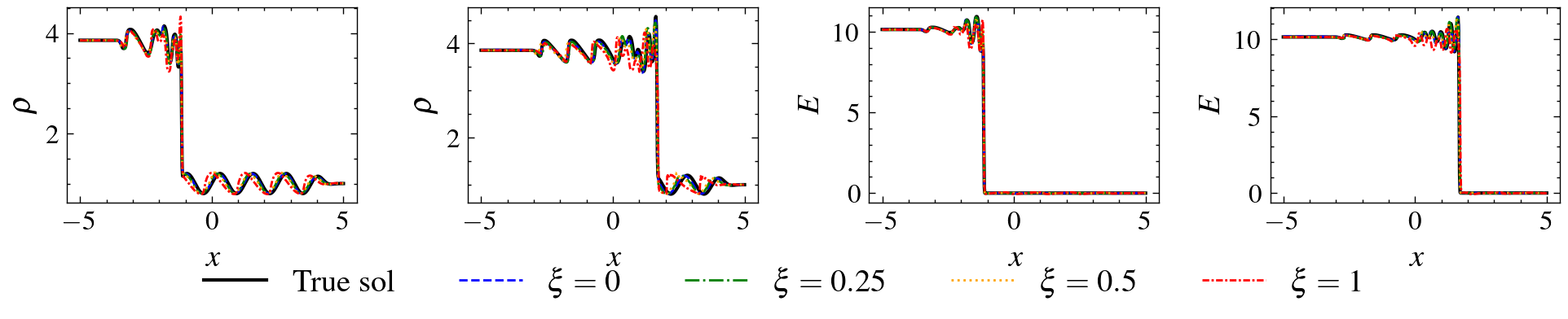}
        \caption{Comparison of the reference solution (black solid line) of density \(\rho\) and energy \(E\) in Euler's equation with the NESCFN model predictions with \(\xi = 0, .25, .5, 1\) in  \eqref{eq:noise_euler}: (left) \(t = .8\) of \(\rho\), (middle-left) \(t = 1.6\) of \(\rho\), (middle-right) \(t =.8\) of \(E\), (right) \(t = 1.6\) of \(E\). }
    \label{fig: Euler Predictions Noise}
\end{figure}

\begin{figure}[h!]
        \centering
        \includegraphics[width=0.6\textwidth]{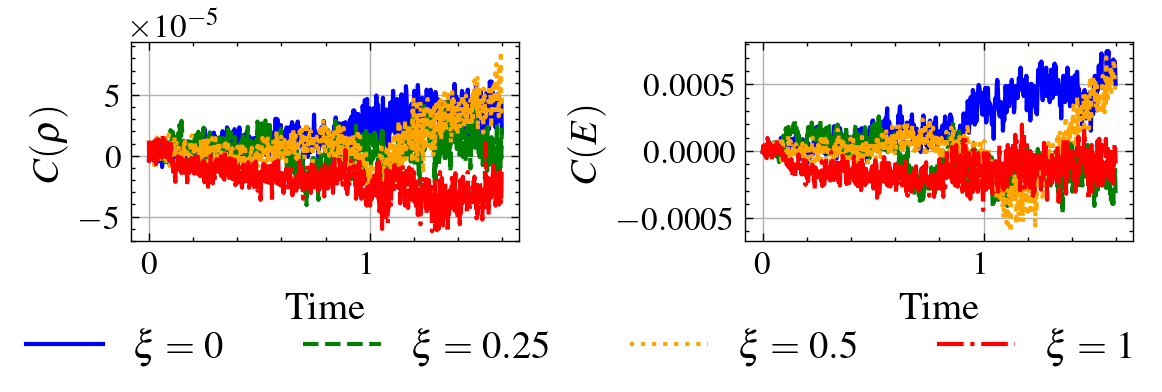}
   \caption{Discrete conserved quantity remainder \(\mathcal{C}(\rho)\)(left), \(\mathcal{C}(E)\)(right), defined in \eqref{eq:conservemetric}, for Euler's equation with noise levels \(\xi = 0, .25, .5, 1\) in  \eqref{eq:noise_euler}.}
        \label{fig: Euler Conservations Noise}
\end{figure}

\begin{figure}[h!]
    \centering
    \begin{subfigure}{0.23\textwidth}
    \centering
        \includegraphics[width=\textwidth]{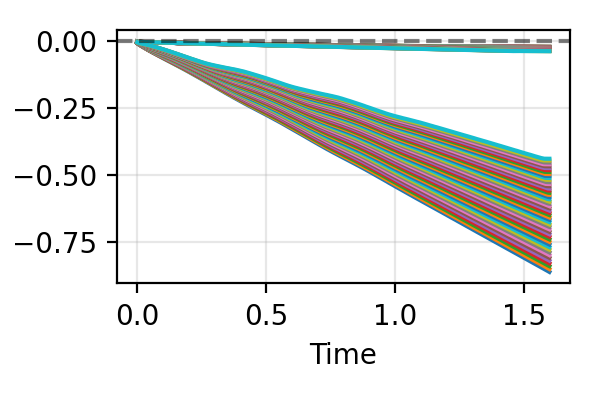}
    \end{subfigure}\hspace{0.1cm}
    \begin{subfigure}{0.23\textwidth}
    \centering
        \includegraphics[width=\textwidth]{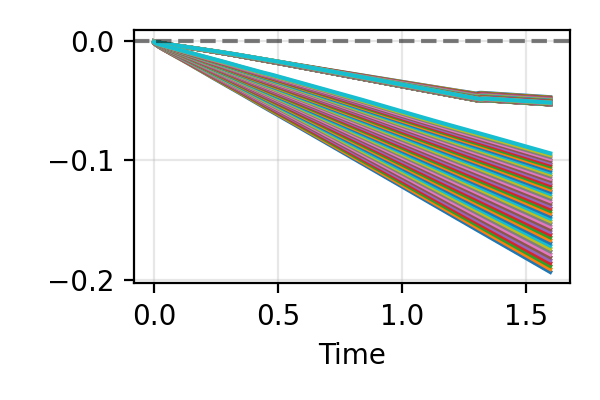}
    \end{subfigure}\hspace{0.1cm}
    \begin{subfigure}{0.23\textwidth}
    \centering
        \includegraphics[width=\textwidth]{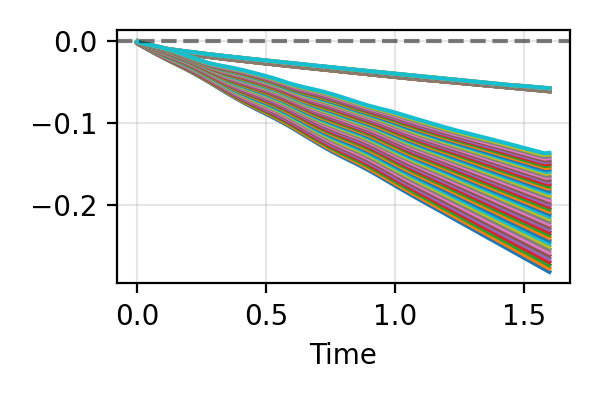}
    \end{subfigure}
    \begin{subfigure}{0.23\textwidth}
    \centering
        \includegraphics[width=\textwidth]{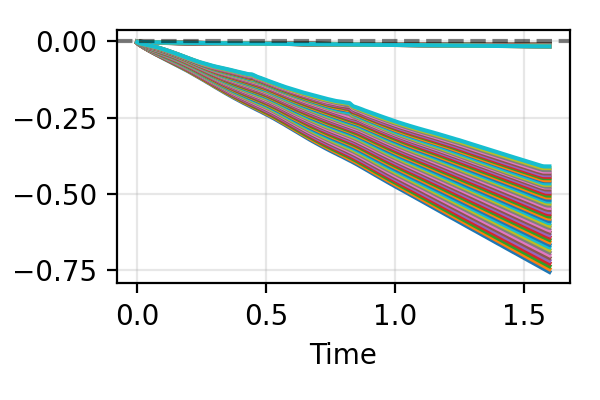}
    \end{subfigure}
    \caption{Discrete entropy remainder $\mathcal{J}([\rho, \rho u, E]^\top)$, defined in \eqref{eq: discrete entropy}, for Euler's equation, evaluated across all testing initial conditions specified in \eqref{eq:entropy_test_1deuler_rod}--\eqref{eq:entropy_test_1d_shuO_osher}. From (left) to (right), the plots correspond to training data with noise levels $\xi = 0, .25, .5$, and $1$, respectively.
    }
    \label{fig: Euler entropy Remainder}
\end{figure}

\Cref{fig: Euler Predictions Noise} shows the model predictions for the density \(\rho\) and the energy \(E\) with different choices of $\xi$ in \eqref{eq:noise_euler}. We omit visualizations of the momentum \(\rho u\) due to its qualitative similarity to \(E\). The prediction accuracy remains comparable to that of the ESCFN model (see \cite[section 5.2.3, Fig. 5.6]{liu2024entropy}), even in very noisy environments. The impact of noise is not negligible however, primarily due to the use of the learned surrogate flux \eqref{eq: neural-entropy stable flux} within the entropy-stable scheme \eqref{eq: entropy stable flux}, which introduces additional diffusion in the high-noise regime. Such behavior is expected -- in the worst case, the learned surrogate approximates the Lax--Friedrichs flux with a constant wave speed, corresponding to the near zero contribution from the network-predicted spectral radius. The robustness and entropy-stability of the NESCFN are still clearly demonstrated, however. The evolution of the discrete conservation remainders \(\mathcal{C}(\rho)\) and \(\mathcal{C}(E)\) shown in \Cref{fig: Euler Conservations Noise}, along with similar behavior for \(\mathcal{C}(\rho u)\), confirms conservation up to numerical tolerance. Additionally, the discrete entropy remainder $\mathcal{J}([\rho, \rho u, E]^\top)$ remains non-positive across all testing initial conditions specified from \eqref{eq:entropy_test_1deuler_rod} to \eqref{eq:entropy_test_1d_shuO_osher}, as illustrated in \Cref{fig: Euler entropy Remainder}.

\subsubsection{2D Burgers' equation}

\begin{figure}[h!]
    \centering
    \includegraphics[width=.8\linewidth]{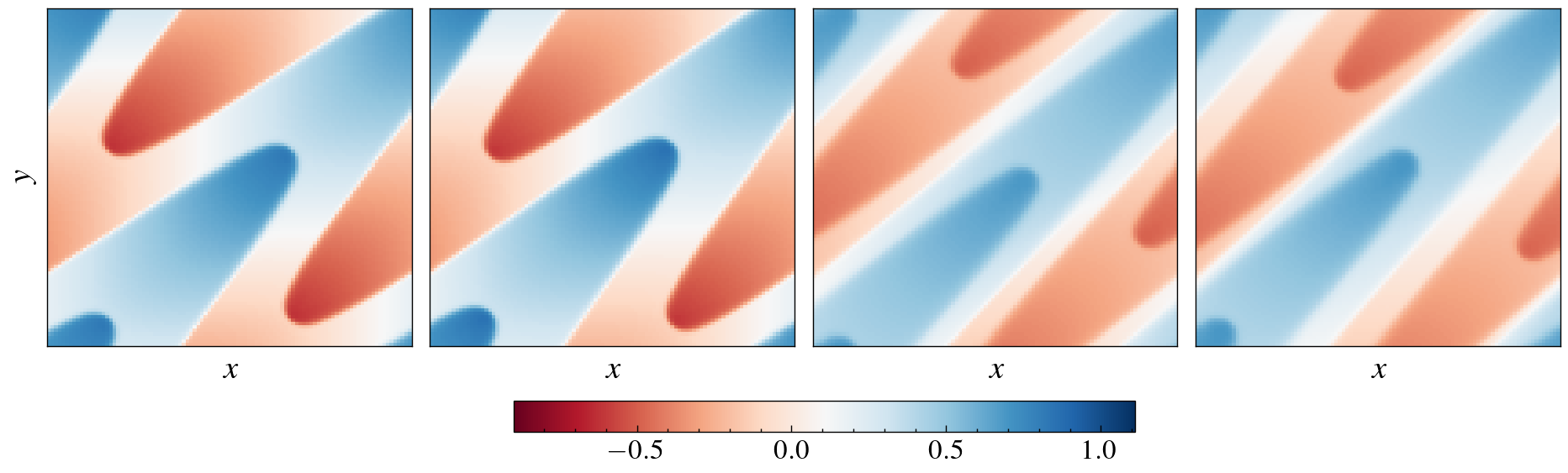}
    \caption{Comparison of the reference solution  of  \(u\) in 2D Burgers' equation \eqref{eq: Burgers' equation 2d} with the NESCFN model predictions for noise coefficient $\xi = 1$: (left) reference solution at t = .8, (middle-left) predictions at t = .8, (middle-right) reference solution at t = 1.6, (right) predictions at t = 1.6.} 
    \label{fig: 2d Burgers Noise}
\end{figure}

We conclude by investigating the impact of noise in the training data for the 2D Burgers' equation \eqref{eq: Burgers' equation 2d}. The noise is introduced following the same procedure as in \eqref{eq: noise burgers}, with noise levels \(\xi = 0, .25, .5, 1\). Recall that $L = L_{\text{train}} = 20$ and $N_{\text{traj}} = 5$. The number of training trajectories as well as the length of training period are smaller due to GPU limitations.  \Cref{fig: 2d Burgers Noise} compares the reference solutions with the predictions of the NESCFN model at times $t = .8$ and $1.6$. The results demonstrate that the NESCFN maintains good predictive accuracy, even with $\xi = 1$.  

\begin{figure}[h!]
    \centering
    \begin{subfigure}{0.23\textwidth}
    \centering
        \includegraphics[width=\textwidth]{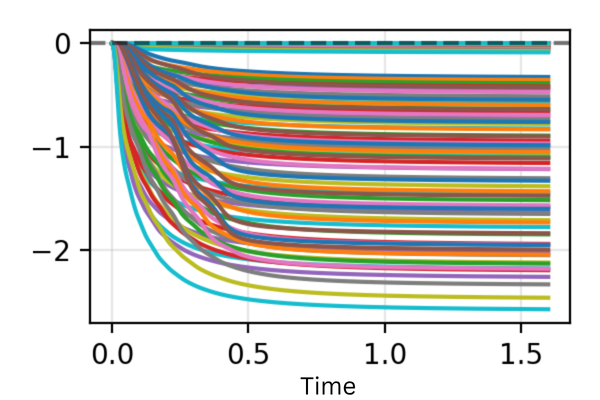}
    \end{subfigure}\hspace{0.1cm}
    \begin{subfigure}{0.23\textwidth}
    \centering
        \includegraphics[width=\textwidth]{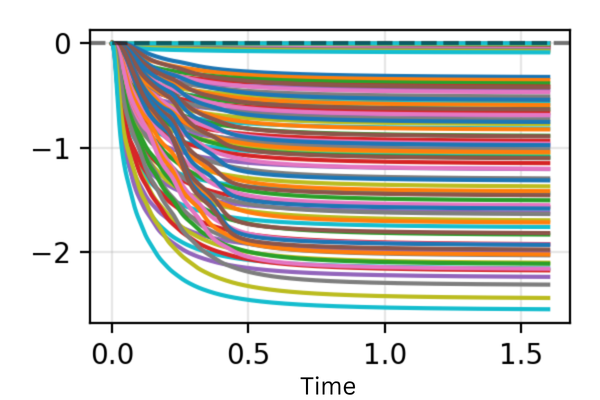}
    \end{subfigure}\hspace{0.1cm}
    \begin{subfigure}{0.23\textwidth}
    \centering
        \includegraphics[width=\textwidth]{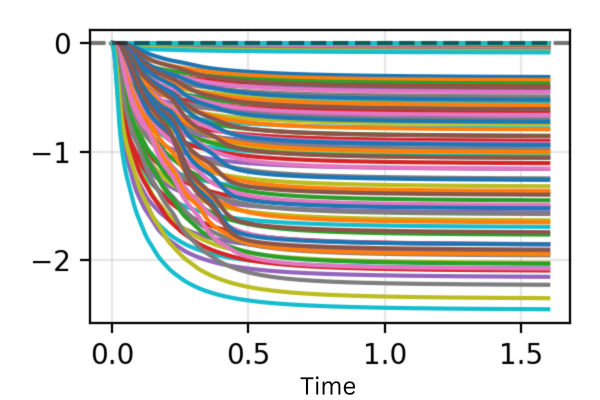}
    \end{subfigure}
    \begin{subfigure}{0.23\textwidth}
    \centering
        \includegraphics[width=\textwidth]{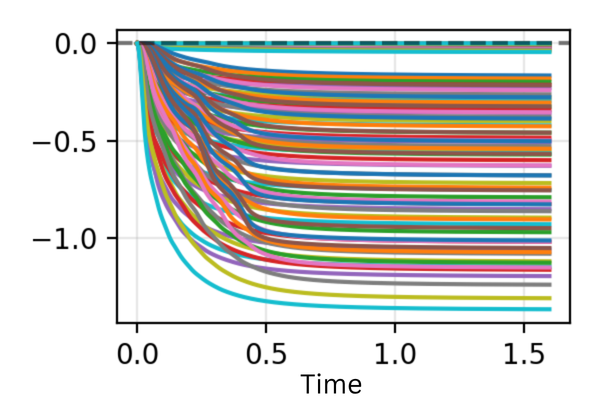}
    \end{subfigure}
    \caption{Discrete entropy remainder $\mathcal{J}(u)$, defined in \eqref{eq: discrete entropy}, for the 2D Burgers' equation, evaluated across all testing initial conditions specified in \eqref{eq:entropy_test_2d_burger1}--\eqref{eq:entropy_test_2dburger3}. From (left) to (right), the plots correspond to training data with noise levels $\xi = 0, .25, .5$, and $1$, respectively.}
    \label{fig: 2d Burgers entropy Remainder}
\end{figure}
\Cref{fig: 2d Burgers entropy Remainder} presents the evolution of the discrete entropy remainder $\mathcal{J}(u)$, defined in \eqref{eq: discrete entropy}. The discrete conserved quantity remainder metric, \(\mathcal{C}(u)\) for the 2D Burgers' equation behaves similarly to the 1D case (see \Cref{fig: 1D Burgers Conservations Noise}), and is therefore omitted. We also observe that applying the neural surrogate to the 2D Burgers' equation introduces more numerical diffusion than expected. This behavior is likely due to difficulties in accurately learning the maximum local wave speed, particularly under noisy and limited training data. Increasing the number of training trajectories and improving data quality may help the model more precisely approximate the spectral radius. A more detailed investigation of these factors is left for future work.

Nonetheless, the results demonstrate that the proposed NESCFN maintains its robustness and entropy-stable properties in the multidimensional setting, even under substantial noise. These findings highlight the generalizability and stability of the NESCFN method as spatial complexity and perturbation levels increase.

\section{Concluding remarks}
\label{sec:conclusion}

This paper introduced the neural {\em entropy-stable} conservative flux form neural network (NESCFN), a data-driven framework for predicting solutions to hyperbolic conservation laws with discontinuities, where both the numerical flux and a convex neural entropy are learned directly from data. By extending the ESCFN proposed in \cite{liu2024entropy}, NESCFN removes the need for predefined entropy pairs and broadens applicability beyond specific numerical schemes. Numerical experiments demonstrate that NESCFN maintains stability and accuracy under noisy training conditions, generalizes to unseen initial conditions, and preserves the entropy inequality throughout long-term predictions. These properties underscore its extrapolatory capability and highlight its potential to bridge entropy-stable theory with real-world data-driven modeling.

Future work may proceed along several important directions. A natural extension is the integration of techniques such as the discontinuous Galerkin method to handle higher-dimensional systems, particularly in the presence of complex geometries, which is an ongoing project. Moreover, while our results demonstrate the empirical effectiveness of the proposed framework, a rigorous theoretical analysis of its convergence and stability remains an important avenue for future investigation.

\section*{Acknowledgements} This work was partially supported by the DOD (ONR MURI) grant \#N00014-20-1-2595,  the DOE ASCR grant \#DE-SC0025555, the NSF grant \#DMS-2513924, and the Ken Kennedy Institute at Rice University
\appendix 

\small
\bibliographystyle{abbrv}
\bibliography{reference}

\begin{thebibliography}{10}

\bibitem{amos2017input}
B.~Amos, L.~Xu, and J.~Z. Kolter.
\newblock Input convex neural networks.
\newblock In {\em International conference on machine learning}, pages 146--155. PMLR, 2017.

\bibitem{Barth1999NumericalMethodsGasDynamic}
T.~J. Barth.
\newblock {\em Numerical Methods for Gasdynamic Systems on Unstructured Meshes}, pages 195--285.
\newblock Springer Berlin Heidelberg, Berlin, Heidelberg, 1999.

\bibitem{iceberg}
G.~R. Bigg, M.~R. Wadley, D.~P. Stevens, and J.~A. Johnson.
\newblock Modelling the dynamics and thermodynamics of icebergs.
\newblock {\em Cold Regions Science and Technology}, 26(2):113--135, 1997.

\bibitem{jax}
J.~Bradbury, R.~Frostig, P.~Hawkins, M.~J. Johnson, C.~Leary, D.~Maclaurin, G.~Necula, A.~Paszke, J.~Vander{P}las, S.~Wanderman-{M}ilne, and Q.~Zhang.
\newblock {JAX}: composable transformations of {P}ython+{N}um{P}y programs, 2018.

\bibitem{chenPDE}
Z.~Chen, V.~Churchill, K.~Wu, and D.~Xiu.
\newblock Deep neural network modeling of unknown partial differential equations in nodal space.
\newblock {\em Journal of Computational Physics}, 449:110782, 2022.

\bibitem{chencfn}
Z.~Chen, A.~Gelb, and Y.~Lee.
\newblock Learning the dynamics for unknown hyperbolic conservation laws using deep neural networks.
\newblock {\em SIAM Journal on Scientific Computing}, 46(2):A825--A850, 2024.

\bibitem{Churchill_2023}
V.~Churchill and D.~Xiu.
\newblock Flow map learning for unknown dynamical systems: Overview, implementation, and benchmarks.
\newblock {\em Journal of Machine Learning for Modeling and Computing}, 4(2):173--201, 2023.

\bibitem{clawpack}
{Clawpack Development Team}.
\newblock Clawpack software, 2020.
\newblock Version 5.7.1.

\bibitem{neuralPDE}
A.~Dulny, A.~Hotho, and A.~Krause.
\newblock Neural{P}{D}{E}: Modelling dynamical systems from data.
\newblock In R.~Bergmann, L.~Malburg, S.~C. Rodermund, and I.~J. Timm, editors, {\em KI 2022: Advances in Artificial Intelligence}, pages 75--89, Cham, 2022. Springer International Publishing.

\bibitem{fjordholm2012arbitrarily}
U.~S. Fjordholm, S.~Mishra, and E.~Tadmor.
\newblock Arbitrarily high-order accurate entropy stable essentially nonoscillatory schemes for systems of conservation laws.
\newblock {\em SIAM Journal on Numerical Analysis}, 50(2):544--573, 2012.

\bibitem{dynonet}
M.~Forgione and D.~Piga.
\newblock \textit{dyno{N}et}: A neural network architecture for learning dynamical systems.
\newblock {\em International Journal of Adaptive Control and Signal Processing}, 35(4):612--626, 2021.

\bibitem{rnn_ode}
K.~Gajamannage, D.~I. Jayathilake, Y.~Park, and E.~M. Bollt.
\newblock {Recurrent neural networks for dynamical systems: Applications to ordinary differential equations, collective motion, and hydrological modeling}.
\newblock {\em Chaos: An Interdisciplinary Journal of Nonlinear Science}, 33(1):013109, 01 2023.

\bibitem{Girard11}
L.~Girard, S.~Bouillon, J.~Weiss, D.~Amitrano, T.~Fichefet, and V.~Legat.
\newblock A new modeling framework for sea-ice mechanics based on elasto-brittle rheology.
\newblock {\em Annals of Glaciology}, 52(57):123–132, 2011.

\bibitem{godlewski2013numerical}
E.~Godlewski and P.-A. Raviart.
\newblock {\em Numerical approximation of hyperbolic systems of conservation laws}, volume 118.
\newblock Springer Science \& Business Media, 2013.

\bibitem{Shu98}
S.~Gottlieb and C.-W. Shu.
\newblock Total variation diminishing {R}unge-{K}utta schemes.
\newblock {\em Math. Comput.}, 67(221):73–85, jan 1998.

\bibitem{Euler}
J.~R. Holton and G.~J. Hakim.
\newblock {\em An Introduction to Dynamic Meteorology}.
\newblock Academic Press, Oxford, 5th edition, 2012.

\bibitem{FarzadAffordableEntropyConsistentEulerFlux}
F.~Ismail and P.~L. Roe.
\newblock Affordable, entropy-consistent {E}uler flux functions ii: Entropy production at shocks.
\newblock {\em Journal of Computational Physics}, 228(15):5410--5436, 2009.

\bibitem{kast2024positional}
M.~Kast and J.~S. Hesthaven.
\newblock Positional embeddings for solving {P}{D}{E}s with evolutional deep neural networks.
\newblock {\em Journal of Computational Physics}, 508:112986, 2024.

\bibitem{pyclaw}
D.~I. Ketcheson, K.~T. Mandli, A.~J. Ahmadia, A.~Alghamdi, M.~{Quezada de Luna}, M.~Parsani, M.~G. Knepley, and M.~Emmett.
\newblock {PyClaw: Accessible, Extensible, Scalable Tools for Wave Propagation Problems}.
\newblock {\em SIAM Journal on Scientific Computing}, 34(4):C210--C231, Nov. 2012.

\bibitem{Adam}
D.~P. Kingma and J.~Ba.
\newblock Adam: {A} method for stochastic optimization.
\newblock In Y.~Bengio and Y.~LeCun, editors, {\em 3rd International Conference on Learning Representations, {ICLR} 2015, San Diego, CA, USA, May 7-9, 2015, Conference Track Proceedings}, 2015.

\bibitem{kurganov2000new}
A.~Kurganov and E.~Tadmor.
\newblock New high-resolution central schemes for nonlinear conservation laws and convection--diffusion equations.
\newblock {\em Journal of computational physics}, 160(1):241--282, 2000.

\bibitem{KTscheme2000}
A.~Kurganov and E.~Tadmor.
\newblock New high-resolution central schemes for nonlinear conservation laws and convection–diffusion equations.
\newblock {\em Journal of Computational Physics}, 160(1):241--282, 2000.

\bibitem{LeFlochEntropyConservativeSchemeArbitrayOrder}
P.~G. LeFloch, J.~M. Mercier, and C.~Rohde.
\newblock Fully discrete, entropy conservative schemes of arbitraryorder.
\newblock {\em SIAM Journal on Numerical Analysis}, 40(5):1968--1992, 2002.

\bibitem{LeVeque02}
R.~J. LeVeque.
\newblock {\em Finite Volume Methods for Hyperbolic Problems}.
\newblock Cambridge Texts in Applied Mathematics. Cambridge University Press, 2002.

\bibitem{deeppropnet}
L.~Liu and W.~Cai.
\newblock Deep{P}rop{N}et -- a recursive deep propagator neural network for learning evolution {P}{D}{E} operators, 2022.

\bibitem{liu2024entropy}
L.~Liu, T.~Li, A.~Gelb, and Y.~Lee.
\newblock Entropy stable conservative flux form neural networks.
\newblock {\em arXiv preprint arXiv:2411.01746}, 2024.

\bibitem{cdeeponet}
L.~Liu, K.~Nath, and W.~Cai.
\newblock A causality-deeponet for causal responses of linear dynamical systems.
\newblock {\em Communications in Computational Physics}, 35(5):1194--1228, 2024.

\bibitem{long2018pdenetlearningpdesdata}
Z.~Long, Y.~Lu, X.~Ma, and B.~Dong.
\newblock {P}{D}{E}-{N}et: Learning {P}{D}{E}s from data, 2018.

\bibitem{Dimitrios2025gorinn}
D.~G. Patsatzis, M.~{di Bernardo}, L.~Russo, and C.~Siettos.
\newblock Gorinns: Godunov-riemann informed neural networks for learning hyperbolic conservation laws.
\newblock {\em Journal of Computational Physics}, 534:114002, 2025.

\bibitem{Elastic-PlasticSeaIce}
R.~S. Pritchard.
\newblock {An Elastic-Plastic Constitutive Law for Sea Ice}.
\newblock {\em Journal of Applied Mechanics}, 42(2):379--384, 06 1975.

\bibitem{roe}
P.~Roe.
\newblock Approximate riemann solvers, parameter vectors, and difference schemes.
\newblock {\em Journal of Computational Physics}, 43(2):357--372, 1981.

\bibitem{SHU1988439}
C.-W. Shu and S.~Osher.
\newblock Efficient implementation of essentially non-oscillatory shock-capturing schemes.
\newblock {\em J. Comput. Phys.}, 77(2):439--471, 1988.

\bibitem{Seung2025tvdnet}
S.~W. Suh, J.~F. MacArt, L.~N. Olson, and J.~B. Freund.
\newblock A {T}{V}{D} neural network closure and application to turbulent combustion.
\newblock {\em Journal of Computational Physics}, 523:113638, 2025.

\bibitem{Tadmor1987EntropyStable}
E.~Tadmor.
\newblock The numerical viscosity of entropy stable schemes for systems of conservation laws. i.
\newblock {\em Mathematics of Computation}, 49:91--103, 1987.

\bibitem{tadmor2016entropy}
E.~Tadmor.
\newblock Entropy stable schemes.
\newblock In {\em Handbook of Numerical Analysis}, volume~17, pages 467--493. Elsevier, 2016.

\bibitem{RoeNet}
Y.~Tong, S.~Xiong, X.~He, S.~Yang, Z.~Wang, R.~Tao, R.~Liu, and B.~Zhu.
\newblock Roenet: Predicting discontinuity of hyperbolic systems from continuous data.
\newblock {\em International Journal for Numerical Methods in Engineering}, 125(6):e7406, 2024.

\bibitem{Shallow_water}
C.~B. Vreugdenhil.
\newblock {\em Numerical Methods for Shallow-Water Flow}.
\newblock Springer, Dordrecht, 1994.

\bibitem{yin2023continuous}
Y.~Yin, M.~Kirchmeyer, J.-Y. Franceschi, A.~Rakotomamonjy, and P.~Gallinari.
\newblock Continuous {P}{D}{E} dynamics forecasting with implicit neural representations.
\newblock In {\em International Conference on Learning Representations (ICLR)}, 2023.

\end{thebibliography}

\end{document}